\documentclass[a4paper, 12pt]{article}
\usepackage{amsfonts}
\usepackage{fancyhdr}
\usepackage{graphicx}
\usepackage{epsfig}
\usepackage{amssymb}
\usepackage{mathrsfs,amsfonts,amsmath}
\usepackage{color}
\usepackage{array}
\usepackage{tabularx,multirow,makecell}
\usepackage{ragged2e}
\usepackage{float}
\usepackage{amsmath,amssymb,amsthm}

\makeatletter

\@addtoreset{equation}{section} \makeatother
\newtheorem{theorem}{Theorem}
\newtheorem{lemma}{Lemma}

\newtheorem{assumption}{Assumption}
\newtheorem{remark}{Remark}
\theoremstyle{definition}
\newtheorem{definition}{Definition}

\DeclareMathOperator{\Law}{Law}
\newcommand{\R}{\mathbb{R}}
\newcommand{\E}{\mathbb{E}}

\newcommand{\Wone}{W_1}
\newcommand{\flow}{\Phi}
\newcommand{\Zt}{\widetilde{Z}}

\newcommand{\Df}{D[\flow]}

\begin{document}

\title{Convergence and exponential stability of splitting Euler method for McKean-Vlasov SDEs
with superlinear drift and H\"older diffusion terms}

\author{Zeping Zhang\quad and\quad Guangqiang Lan\footnote{Corresponding author: Email: langq@buct.edu.cn.}
\\ \scriptsize College of Mathematics and Physics, Beijing University of Chemical Technology, Beijing, China}
\date{}
\maketitle

\begin{abstract}
We consider one-dimensional McKean-Vlasov stochastic differential
equations (MVSDEs) whose drift terms are superlinear and diffusion term are $\beta+\frac{1}{2}$ H\"{o}lder continuous ($\beta\in[0,\frac{1}{2}]$). By decomposing $f(x,\mu)=f_1(x,\mu)+f_2(x)$,
where $f_1$ is globally Lipschitz continuous with respect to the state and measure
variables and $f_2$ is a superlinear, a splitting Euler method is introduced to discretize the given equation. $L^1$ convergence of splitting Euler method to the corresponding interacting particle system are obtained with the convergence rates $\Delta^\beta+N^{-1/2}$ for $\beta\in(0,\frac12]$ and $\log(1/\Delta)^{-1}+N^{-1/2}$ for $\beta=0$. Moreover, $L^p$ uniform convergence are also obtained with lower rates in different cases.  Under additional strict dissipativity condition on $f_2$, we also derive $L^1$ exponential decay estimate for the exact solution and $L^2$ exponential decay estimation of splitting Euler method with the exponentially
tamed-diffusion scheme.
\end{abstract}



\section{Introduction}

MVSDEs, in which
the coefficients depend on the law of the solution, arise as
mean-field limits of weakly interacting particle systems and appear
in kinetic theory, mean-field games, and financial modelling; see
e.g.\ \cite{CD,DRS} and the references therein. In this paper, we consider the following
one-dimensional MVSDE:
\begin{equation}\label{eq:MVSDE}
    dX_t = f(X_t,\mu_t)\,dt + g(X_t)\,dW_t,\qquad t\ge 0,
    \qquad \mu_t=\Law(X_t),
\end{equation}
whose drift can be decomposed as
\begin{equation}\label{eq:decomp}
    f(x,\mu)=f_1(x,\mu)+f_2(x),
\end{equation}
where $f_1$ carries the measure dependence and is globally
Lipschitz, while $f_2$ is a state-dependent drift with superlinear
(polynomial) growth, satisfying a one-sided Lipschitz and a
dissipativity condition. The diffusion $g$ is only
H\"older continuous. A prototype is
$f_2(x)=-x^{2k+1}$ combined with a square-root-type diffusion
$g(x)=|x|^{1/2}$.

Two difficulties are typical for this kind of equations. Firstly, explicit Euler
approximations of SDEs with superlinear coefficients may diverge, as pointed out in
\cite{HJKdiv}; established alternatives include taming
\cite{HJK,Kieu2019,Kieu2022}, truncation \cite{Mao}, and adaptive
time stepping \cite{FG,Kieu2022}. Secondly, since diffusion coefficients are only
H\"older continuous, the Yamada--Watanabe analysis
for Euler approximations gives rates of order $\Delta^\beta$ for
$\beta>0$ and logarithmic order at the endpoint $\beta=0$, under the
conditions of \cite{GR}. Tamed Euler methods for
McKean--Vlasov equations with superlinear drift and H\"older diffusion
are considered in \cite{LiuShiWu}. For long-time behavior, \cite{LiuGaoGuo} study moment stability, infinite-horizon propagation of chaos,
and mean-square stability of a modified tamed Euler scheme under H\"older diffusion.
More broadly, explicit and implicit Euler-type schemes for McKean--Vlasov equations
with superlinear drift have been analysed in \cite{DRS,ChenReisSSM,ChenReisSSM2024}. 
Euler method for  McKean--Vlasov equations with non-Lipschitz coefficients is considered with linear growth condition in \cite{DQ}. 
Split-step methods have also been developed directly for
McKean--Vlasov equations with superlinear drift and Lipschitz diffusion
\cite{ChenReisSSM,ChenReisSSM2024}; explicit and implicit particle-based
Euler methods for superlinear McKean--Vlasov equations are studied in
\cite{DRS} for Lipschitz continuous (both in space and in measure) diffusion case.

The present method differs from these schemes in using the exact ODE flow
of the dissipative drift component and in treating the H\"older diffusion through a Yamada--Watanabe error analysis.
The main numerical device in this paper is a flow-based splitting.
Splitting idea is classical in deterministic numerical analysis \cite{HLW} (see also \cite{B})
and has been studied for SDEs with locally Lipschitz, one-sided
Lipschitz drifts and additive noise in \cite{BSTT}. In \cite{BSTT}, the authors
establish mean-square convergence of the splitting method with order one for
additive-noise SDEs and preserves structural properties by using the fundamental theorem of \cite{TZ}.

Inspired by \cite{BSTT}, we will introduce the splitting Euler method of MVSDE \eqref{eq:MVSDE} with H\"older continuous coefficient.
Our scheme integrates the
selected dissipative component $f_2$ through its exact ODE flow and
handles a merely H\"older-continuous diffusion by a flow-corrected
It\^o interpolation and a Yamada--Watanabe error argument.

Our main results are as follows.
\begin{enumerate}
\item \textbf{Flow-corrected interpolation.} The key device (Section \ref{sec:scheme}) is a
    continuous interpolation $\Zt_t=\flow_{t-t_k}(Z_t^{\mathrm{EM}})$
    which is a It\^o process whose drift contains
    $f_2(\Zt_t)$ itself. The one-sided Lipschitz estimate for $f_2$
    therefore applies directly in the Yamada--Watanabe error
    equation, and no truncation error, stopping time, or truncation-tail estimate is
    needed.
    \item \textbf{Finite-time error estimation} (Theorem \ref{thm:main}):
    For the non-interacting particle system $X_t^i$ and the proposed splitting Euler scheme
    $\tilde{Z}_t^i$, without using the method of propagation of chaos, we establish
    $L^1$ error bounds $\sup_{1\le i\le N}\sup_{t\in[0,T]}\E|X_t^i-\tilde{Z}_t^i|$ of order
    $\Delta^\beta+N^{-1/2}$ for $0<\beta\le\frac12$ and
    $1/\log(1/\Delta)+N^{-1/2}$ for $\beta=0$.
    \item \textbf{Uniform error estimation} (Theorem \ref{thm:pathsup1}):
    We also establish $L^p$ uniform error bounds $\sup_{1\le i\le N}\E\sup_{t\in[0,T]}|X_t^i-\tilde{Z}_t^i|^p$
    under the same assumptions with relatively lower rate.
    Moreover, our $L^1$ and $L^2$ rates are better than that of tamed Euler method in \cite{LiuShiWu} in some special cases (see Remark \ref{zhu} for interpretations).
        \item \textbf{Long-time second-moment decay} (Theorem \ref{thm:stab}):
    under equilibrium and strict dissipativity assumptions, the exponentially
    tamed scheme satisfies an explicit exponential decay estimate for its
    second moment. 
\end{enumerate}

The paper is organized as follows. Section \ref{sec:model} collects
the assumptions and preliminaries, and the exact-flow estimates are established. Then in Section \ref{sec:scheme}, the
splitting scheme and its flow-corrected interpolation are defined, the moment bounds are also obtained. In Section
\ref{sec:conv}, main convergence theorems are proved by using Yamada-Watanabe approximation. Section \ref{sec:stab} studies long-time exponential
decay of the exact and numerical solutions. Finally, we give a concluding remark in Section \ref{jiesu}.

\section{Model, assumptions and preliminaries}\label{sec:model}

Throughout, $\mathcal{P}_1(\R)$ denotes the space of probability
measures on $\R$ with finite first moment, endowed with the
Wasserstein-$p$ ($p\ge1$) distance
\[
    W_p^p(\mu,\nu)=\inf_{\pi\in\Pi(\mu,\nu)}
    \int_{\R\times\R}|x-y|^p\,\pi(dx,dy),
    \qquad
    \|\mu\|_p^p=\int_\R |x|^p\,\mu(dx)
\]
where $\Pi(\mu,\nu)$ is the set of all of the coupling probability measures of  $\mu$ and $\nu$.

\begin{assumption}\label{ass:A1}
There exists $L>0$ such that for every $x,y\in\R$ and
$\mu,\nu\in\mathcal{P}_1(\R)$,
\[
    |f_1(x,\mu)-f_1(y,\nu)|
    \le L\bigl(|x-y|+\Wone(\mu,\nu)\bigr).
\]
\end{assumption}
Obviously, the above Lipschitz condition implies
$$|f_1(x,\mu)|\le L\bigl(1+|x|+\|\mu\|_1\bigr).$$

\begin{assumption}\label{ass:A2}
$f_2\in C^2(\R)$, and there exist constants $L>0$, $l\ge0$,
$\lambda>0$ and $C>0$ such that for every
$x,y\in\R$,
\[
\begin{aligned}
    |f_2(x)-f_2(y)|&\le L\bigl(1+|x|^l+|y|^l\bigr)|x-y|,\\
    |f_2'(x)|&\le C\bigl(1+|x|^l\bigr),\\
    |f_2''(x)|&\le C\bigl(1+|x|^{2l}\bigr).
\end{aligned}
\]
\[
    (x-y)\bigl(f_2(x)-f_2(y)\bigr)\le \lambda|x-y|^2.
\]
\end{assumption}

\begin{assumption}\label{ass:A3}
There exist $\rho>0$ and $\beta\in[0,\frac12]$ such that for
every $x,y\in\R$,
\[
    |g(x)-g(y)|\le \rho|x-y|^{1/2+\beta}.
\]
\end{assumption}


 According to \cite{LiuShiWu}, equation
\eqref{eq:MVSDE} admits a unique strong solution under Assumptions \ref{ass:A1}-\ref{ass:A3} with
$\E[\sup_{0\le t\le T}|X_t|^p]\le C_{T,p}$, if $\E|X_0|^{p(l+1)}<\infty$.

Let $(W^i)_{i\ge1}$ be independent standard Brownian motions and
$(X_0^i)_{i\ge1}$ i.i.d.\ copies of $X_0$. The corresponding non-interacting and interacting particle systems of \eqref{eq:MVSDE} are
\begin{equation}\label{eq:NIPS}
    dX^i_t = f(X^i_t,\mu_t)\,dt + g(X^i_t)\,dW^i_t,
    \qquad \mu_t=\Law(X^i_t)=\Law(X_t),
\end{equation}
and
\begin{equation}\label{eq:IPS}
    dX_t^{i,N}=f\bigl(X_t^{i,N},\mu_t^{X,N}\bigr)dt
    +g\bigl(X_t^{i,N}\bigr)dW_t^i,
    \qquad
    \mu_t^{X,N}=\frac1N\sum_{j=1}^N\delta_{X_t^{j,N}}.
\end{equation}

By Theorem 3.1 in \cite{LiuShiWu}, the following empirical-measure estimate is the standard one-dimensional
$W_1$ empirical bound; see Fournier--Guillin \cite{FGuillin}.

\begin{lemma}\label{lem:empirical}
Suppose Assumptions \ref{ass:A1}-\ref{ass:A3} hold. If $\E|X_0|^{q(l+1)}<\infty$, for every $T>0$,
\[
    \sup_{0\le s\le T}
    \E W_p^p\bigl(\mu_s,\bar\mu_s^{X,N}\bigr)\le C_T(N^{-1/2}+N^{-\frac{q-p}{q}}), q\neq 2p,
    \quad
    \bar\mu_s^{X,N}=\frac1N\sum_{j=1}^N\delta_{X_s^j}.
\]
Specially, for $\E|X_0|^{2(l+1)}<\infty$,
\[
    \sup_{0\le s\le T}
    \E\Wone\bigl(\mu_s,\bar\mu_s^{X,N}\bigr)\le C_TN^{-1/2}.
\]
\end{lemma}

Consider the ordinary differential equation
\begin{equation}\label{eq:ode}
    \frac{dy(t)}{dt}=f_2(y(t)),\qquad y(0)=x.
\end{equation}
Denote $\flow_t(x)$ be the solution flow to (\ref{eq:ode}) at $t$ with initial value $x$. That is
\begin{equation}\flow_t(x)=y(t,x).\end{equation}

We have the following Lemma.
\begin{lemma}\label{lem:flow}
Let Assumption \ref{ass:A2} hold. Then there is a unique solution flow $\flow_t(x)$ to (\ref{eq:ode}). Moreover, for every $x,y\in\R$ and
$t\in[0,1]$:
\begin{enumerate}
    \item[(F1)]     $|\flow_t(x)-\flow_t(y)|\le e^{\lambda t}|x-y|$;
    \item[(F2)]    $|\flow_t(x)-x|\le Ct\bigl(1+|x|^{l+1}\bigr)$;
    \item[(F3)]   Denote $ R(x,t):=\flow_t(x)-x-f_2(x)t.$ Then
    \[|R(x,t)|\le Ct^2\bigl(1+|x|^{2l+1}\bigr);
    \]
    \item[(F4)] Denote $D[\flow]_t(x):=\frac{\partial}{\partial x}\flow_t(x)$ and $D^2[\flow]_t(x):=\frac{\partial^2}{\partial x^2}\flow_t(x)$. Then
    \[|D[\flow]_t(x)|\le e^{\lambda t}, \quad  |D^2[\flow]_t(x)|\le Ct\bigl(1+|x|^{2l}\bigr),\]
       \[ |D[\flow]_t(x)-1|\le Ct\bigl(1+|x|^l\bigr)
      \]
    and $$D[\flow]_t(x)f_2(x)=f_2(\flow_t(x)).$$
  \item[(F5)] For every $p\ge2$ there
    are constants $c_p,C_p>0$ such that
    \[
      |\flow_t(x)|^p\le e^{c_p t}|x|^p+C_p t,
      \qquad 0\le t\le1.
    \]
\end{enumerate}
\end{lemma}

\begin{proof}
Since $f_2$ is locally Lipschitz, the ODE has a unique local solution. The one-sided Lipschitz and dissipativity
bounds below prevent finite-time explosion, so this solution flow $\flow_t(x)$ is well defined for all $t\ge0$.

\smallskip
\noindent\textbf{(F1)}
Let $u(t)=\flow_t(x)-\flow_t(y)$. Then
\[
 \frac{d}{dt}|u(t)|^2
 =2u(t)\bigl(f_2(\flow_t(x))-f_2(\flow_t(y))\bigr)
 \le 2\lambda|u(t)|^2.
\]
Gronwall's inequality gives $|u(t)|^2\le e^{2\lambda t}|x-y|^2$,
taking square roots proves (F1).

\smallskip
\noindent\textbf{(F2)} The integral equation for the flow gives
\[
 |\flow_t(x)-x|\le\int_0^t|f_2(\flow_s(x))|\,ds.
\]
By the polynomial local-Lipschitz bound in Assumption~\ref{ass:A2},
with $y=0$,
\[
 |f_2(z)|\le |f_2(0)|+L(1+|z|^l)|z|.
\]
On the other hand, The one-sided Lipschitz condition yields
\[
 \flow_t(x)\bigl(f_2(\flow_t(x))-f_2(0)\bigr)
 \le \lambda|\flow_t(x)|^2.
\]
Consequently, using $2ab\le a^2+b^2$,
\[
 \begin{aligned}
 \frac{d}{dt}|\flow_t(x)|^2
 &=2\flow_t(x)f_2(\flow_t(x))\\
 &\le 2\lambda|\flow_t(x)|^2
       +2|f_2(0)|\,|\flow_t(x)|\\
 &\le (2\lambda+1)|\flow_t(x)|^2+|f_2(0)|^2.
 \end{aligned}
\]
Set $a=2\lambda+1\ge1$. Applying Gronwall's inequality in its
inhomogeneous form gives
\[
 |\flow_t(x)|^2\le e^{at}|x|^2
       +\frac{e^{at}-1}{a}|f_2(0)|^2
 \le e^{at}\bigl(|x|^2+|f_2(0)|^2\bigr).
\]
Thus, for any $s\le t\le 1$, $|\flow_s(x)|\le C(1+|x|)$ where $C=e^{\frac{a}{2}}(1\vee|f_2(0)|^2)$.

Then
$|f_2(y_x(s))|\le C(1+|x|^{l+1})$ for $s\in[0,1]$, and integration over $[0,t]$
proves (F2).

\smallskip
\noindent\textbf{(F3)}
By the integral equation,
\[
 R(x,t)=\int_0^t\bigl(f_2(\flow_s(x))-f_2(x)\bigr)\,ds.
\]
The local-Lipschitz bound in Assumption~\ref{ass:A2} therefore implies
\[
 |R(x,t)|\le L\int_0^t
 (1+|\flow_s(x)|^l+|x|^l)|\flow_s(x)-x|\,ds.
\]
For $s\le t\le1$, (F2) gives $|\flow_s(x)|\le C(1+|x|)$ and
$|\flow_s(x)-x|\le Cs(1+|x|^{l+1})$. Hence the integrand is bounded by
$Cs(1+|x|^{2l+1})$. Integrating $s$ from $0$ to $t$ proves
$|R(x,t)|\le Ct^2(1+|x|^{2l+1})$.

\smallskip
\noindent\textbf{(F4)}
Because $f_2\in C^2$, differentiation of the ODE with respect to its
initial value is justified and $J_1(t):=D[\flow]_t(x)$ satisfies
\[
 J_1'(t)=f_2'(\flow_t(x))J_1(t),\qquad J_1(0)=1.
\]
The one-sided Lipschitz condition implies $f_2'(z)\le\lambda$ for
all $z$ (take the limit $y\to x$ in that condition). Thus
\[
 J_1(t)=\exp\!\left(\int_0^t f_2'(\flow_s(x))\,ds\right),
 \qquad |J_1(t)|\le e^{\lambda t}.
\]
Next, since $J_1(0)=1$,
\[
 |J_1(t)-1|
 \le\int_0^t |f_2'(\flow_s(x))|\,|J_1(s)|\,ds.
\]
The derivative-growth assumption and (F2) imply
$|f_2'(\flow_s(x))|\le C(1+|x|^l)$ for $s\in[0,1]$; also
$|J_1(s)|\le e^{\lambda}$. It follows that
$|D[\flow]_t(x)-1|\le Ct(1+|x|^l)$.

For the second derivative, set $J_2(t):=D^2[\flow]_t(x)$. A second
differentiation of the flow equation gives
\[
 J_2'(t)=f_2''(\flow_t(x))J_1(t)^2+f_2'(\flow_t(x))J_2(t),
 \qquad J_2(0)=0.
\]
Variation of constants, together with the formula for $J_1$, yields
\[
 J_2(t)=J_1(t)\int_0^t f_2''(\flow_s(x))J_1(s)\,ds.
\]
Using $|f_2''(z)|\le C(1+|z|^{2l})$, (F2), and the bound on $J_1$,
we obtain
\[
 |D^2[\flow]_t(x)|\le Ct(1+|x|^{2l}),\qquad 0\le t\le1.
\]

Finally, the flow property
$\flow_t(\flow_s(x))=\flow_{t+s}(x)$, differentiated with respect
to $s$ at $s=0$, gives
$D[\flow]_t(x)f_2(x)=f_2(\flow_t(x))$.

\noindent\textbf{(F5)}
For $p\ge2$, the chain rule and the one-sided Lipschitz condition in
Assumption~\ref{ass:A2} give
\[
 \begin{aligned}
 \frac{d}{dt}|\flow_t(x)|^p
 &=p|\flow_t(x)|^{p-2}\flow_t(x)f_2(\flow_t(x))\\
 &\le p\lambda|\flow_t(x)|^p
       +p|\flow_t(x)|^{p-1}|f_2(0)|.
 \end{aligned}
\]

Therefore, Young's inequality implies
\[
 \frac{d}{dt}|\flow_t(x)|^p\le (p\lambda+p-1)|\flow_t(x)|^p+|f_2(0)|^p.
\]
Then by Gronwall's inequality and mean value theorem, it follows that
\[
 \aligned |\flow_t(x)|^p&\le e^{(p\lambda+p-1)t}|x|^p
       +\frac{|f_2(0)|^p}{p\lambda+p-1}(e^{(p\lambda+p-1)t}-1)\\&
       =e^{(p\lambda+p-1)t}|x|^p
       +\frac{|f_2(0)|^p}{p\lambda+p-1}e^{(p\lambda+p-1)\xi}t\\&
       \le e^{(p\lambda+p-1)t}|x|^p
       +\frac{|f_2(0)|^p e^{p\lambda+p-1}}{p\lambda+p-1}t.\endaligned
\]
\end{proof}

\section{Splitting Euler scheme and moment bounds}\label{sec:scheme}

\begin{definition}\label{def:scheme}
Let $\Delta\in(0,1)$, $t_k=k\Delta$, and
$\hat\mu_k^Z=\frac1N\sum_{j=1}^N\delta_{Z_k^j}$. Set
$Z_0^i=X_0^i$ and, for $k\ge0$,
\begin{equation}\label{eq:scheme}
    Z_{k+1}^i
    =\flow_\Delta\Bigl(
    Z_k^i+f_1\bigl(Z_k^i,\hat\mu_k^Z\bigr)\Delta
    +g\bigl(Z_k^i\bigr)\,\Delta W_k^i
    \Bigr),
    \qquad
    \Delta W_k^i=W_{t_{k+1}}^i-W_{t_k}^i.
\end{equation}
\end{definition}
The central technique of this paper is the following
continuous interpolation. For $t\in[t_k,t_{k+1})$ define the
Euler-stage process
\begin{equation}\label{eq:zem}
    Z^{i,\mathrm{EM}}_t=Z^i_k+f_1\bigl(Z^i_k,\hat\mu_k^Z\bigr)(t-t_k)
    +g\bigl(Z^i_k\bigr)\bigl(W^i_t-W^i_{t_k}\bigr),
\end{equation}
and the \emph{flow-corrected interpolation}
\begin{equation}\label{eq:interp}
    \Zt_t=\flow_{t-t_k}\bigl(Z^{i,\mathrm{EM}}_t\bigr),
    \qquad t\in[t_k,t_{k+1})
\end{equation}
is called \textbf{splitting Euler scheme}.

\begin{lemma}\label{lem:interp}
The process $\Zt^i_t$ is continuous on $[0,T]$ and coincides with the
scheme (\ref{eq:zem}) at the grid points, i.e., $\Zt^i_{t_k}=Z^{i,\mathrm{EM}}_{t_k}=Z^i_k$, and satisfies on each
interval $(t_k,t_{k+1})$
\begin{equation}\label{eq:interp-sde}
\begin{aligned}
    d\Zt^i_t
    ={}&\Bigl[f_2(\Zt^i_t)
    +D[\flow]_{t-t_k}(Z^{i,\mathrm{EM}}_t)\,f_1\bigl(Z^i_k,\hat\mu_k^Z\bigr)\\
    &\qquad +\frac12D^2[\flow]_{t-t_k}(Z^{i,\mathrm{EM}}_t)\,
    g\bigl(Z^i_k\bigr)^2\Bigr]dt\\
    &\quad +D[\flow]_{t-t_k}(Z^{i,\mathrm{EM}}_t)\,g\bigl(Z^i_k\bigr)\,dW^i_t.
\end{aligned}
\end{equation}
\end{lemma}

\begin{proof}
Continuity: at $t=t_{k+1}$,
$\Zt^i_{t_{k+1}^-}=\flow_\Delta(Z^{i,\mathrm{EM}}_{t_{k+1}^-})=Z^i_{k+1}$ by
\eqref{eq:scheme}, while
$\Zt^i_{t_{k+1}}=\flow_0(Z^{i,\mathrm{EM}}_{t_{k+1}})=Z^i_{k+1}$. The dynamics
\eqref{eq:interp-sde} can be achieved by It\^o's formula applied to
$\flow_{t-t_k}(Z^{i,\mathrm{EM}}_t)$, using
$\frac{\partial}{\partial t}\flow_t(x)=f_2(\flow_t(x))$.
\end{proof}

\begin{lemma}\label{lem:moments}
Let Assumptions \ref{ass:A1}-\ref{ass:A3} hold. If $\E |X_0|^{p}<\infty,$ then for every
$T>0$ and every $p>0$ there exists $C_{T,p}$, independent of
$\Delta$ and $N$, such that
\[
    \sup_{0\le k\le [T/\Delta]}\sup_{0< \Delta\le \Delta^*}\sup_{1\le i\le N}\E|Z_k^i|^p
    +\sup_{0< \Delta\le \Delta^*}\sup_{1\le i\le N}\sup_{0\le t\le T}\E|\Zt_t^i|^p
    \le C_{T,p}.
\]
\end{lemma}

\begin{proof}
Fix $p\ge2$ and set
$a_k^{(p)}=\sup_{0\le j\le k}\sup_{1\le i\le N}\E|Z_j^i|^p$.
Write the Euler stage as
\[
Y_k^i=Z_k^i+f_1(Z_k^i,\hat\mu_k^Z)\Delta+g(Z_k^i)\Delta W_k^i.
\]
The linear-growth assumptions on $f_1,g$ imply
\[
\E\bigl[|Y_k^i|^p\bigr]
\le (1+C_p(\Delta^p+\Delta^\frac{p}{2}))\E|Z_k^i|^p
 +C_p\Delta^p\bigl(1+\E\|\hat\mu_k^Z\|_1^p\bigr).
\]
Taking expectations and using Jensen's inequality,
$\E\|\hat\mu_k^Z\|_1^p\le a_k^{(p)}$, gives
\[
\sup_{0\le j\le k}\sup_{1\le i\le N}\E|Y_j^i|^p\le(1+C_p\Delta^\frac{p}{2})a_k^{(p)}+C_p\Delta^p.
\]
Now applying (F5) to $Y_k^i$ yields:
\[
\E|Z_{k+1}^i|^p=\E|\flow_\Delta(Y_k^i)|^p
\le e^{(p\lambda+p-1)\Delta}\E|Y_k^i|^p
       +\frac{|f_2(0)|^p e^{p\lambda+p-1}}{p\lambda+p-1}\Delta.
\]
Consequently,
\[
\aligned a_{k+1}^{(p)}&\le(1+C_p\Delta)e^{(p\lambda+p-1)\Delta}a_k^{(p)}+(C_pe^{p\lambda+p-1}+\frac{|f_2(0)|^p e^{p\lambda+p-1}}{p\lambda+p-1})\Delta\\&
\le e^{c\Delta}a_k^{(p)}+C\Delta\endaligned
\]
where $c=C_p+p\lambda+p-1, C=C_pe^{p\lambda+p-1}+\frac{|f_2(0)|^p e^{p\lambda+p-1}}{p\lambda+p-1}$.

Therefore,
$$\aligned a_{k}^{(p)}&\le e^{ck\Delta}a_0^{(p)}+C\Delta(1+e^{c\Delta}+\cdots+e^{c(k-1)\Delta})\\&
\le e^{ck\Delta}a_0^{(p)}+C\Delta\frac{e^{ck\Delta}}{e^{c\Delta}-1}\\&
\le e^{cT}(a_0^{(p)}+\frac{C}{c}).\endaligned$$
Since the constant on the right hand side is independent of $\Delta$ and $k,$
$$\sup_{0< \Delta\le \Delta^*}\sup_{0\le k\le [T/\Delta]}\sup_{1\le i\le N}\E|Z_k^i|^p
        \le C_{T,p}$$ holds for $C_{T,p}=e^{cT}(a_0^{(p)}+\frac{C}{c}).$

Same moment estimate procedure with $t-t_k\le\Delta$ gives
$\E|Z^{i,\mathrm{EM}}_t|^p\le C_{T,p}$. Applying (F5) again to
$\Zt_t^i=\flow_{t-t_k}(Z^{i,\mathrm{EM}}_t)$ gives
$$\sup_{0< \Delta\le \Delta^*}\sup_{1\le i\le N}\sup_{0\le t\le T}\E|\Zt_t^i|^p\le C_{T,p}.$$

For $0<p<2$, the same conclusion follows by Lyapunov's inequality.
\end{proof}

\begin{lemma}\label{lem:increments}
Under the assumptions of Lemma \ref{lem:moments}, if $\E |X_0|^{p(l+1)}<\infty,$ then for any
$t\in[t_k,t_{k+1})$ and $p>0$,
\[
    \sup_{1\le i\le N}\E|Z^{i,\mathrm{EM}}_t-Z^i_k|^p\le C_p\Delta^{p/2},
    \qquad
    \sup_{1\le i\le N}\E|\Zt^i_t-Z^i_k|^p\le C_p\Delta^{p/2}.
\]
\end{lemma}

\begin{proof}
The first bound follows from \eqref{eq:zem}, the linear growth of
$f_1$ and $g$, the moment estimates in Lemma \ref{lem:moments}, and
$\E|W^i_t-W^i_{t_k}|^p\le C_p\Delta^{p/2}$.

For the second, write
\[
    |\Zt^i_t-Z^i_k|
    \le\bigl|\flow_{t-t_k}(Z^{i,\mathrm{EM}}_t)-Z^{i,\mathrm{EM}}_t\bigr|+|Z^{i,\mathrm{EM}}_t-Z_k|.
\]
By (F2), the first term is bounded by
$C\Delta(1+|Z^{i,\mathrm{EM}}_t|^{l+1})$. Then Lemma
\ref{lem:moments} implies its $p$-th moment is bounded by
$C_{T,p}\Delta^p$. For $0<p\le1$
use $|a+b|^p\le|a|^p+|b|^p$; for $p>1$ use the standard inequality
$|a+b|^p\le2^{p-1}(|a|^p+|b|^p)$.

So we have
$$ \E|\Zt^i_t-Z^i_k|^p\le C_{T,p}\Delta^p+ C_p\Delta^{p/2}$$

Since $0<\Delta<1$,
$\Delta^p\le\Delta^{p/2}$. We complete the proof.
\end{proof}

\section{Convergence analysis}\label{sec:conv}

Write
$\kappa(t):=t_k$ for $t\in[t_k,t_{k+1})$ and, on $[t_k,t_{k+1})$,
\begin{equation}\label{eq:Abar}
\bar A^i(t):=f_1(Z^i_{\kappa(t)},\hat\mu_{\kappa(t)}^Z),\qquad
\bar B^i(t):=g(Z^i_{\kappa(t)}),
\end{equation}
so that
$$Z^{i,\mathrm{EM}}_t=Z^i_{\kappa(t)}+\bar A^i(t)(t-\kappa(t))+\bar B^i(t)(W^i_t-W^i_{\kappa(t)})$$
and $\Zt^i_t=\flow_{t-\kappa(t)}(Z^{i,\mathrm{EM}}_t)$.

We are now in the position to state the main convergence theorem.

\begin{theorem}\label{thm:main}
Let Assumptions \ref{ass:A1}-\ref{ass:A3} hold with $0\le\beta\le\tfrac12$ and let $0<\Delta<1$. If $\E |X_0|^{2(l+1)}<\infty,$ then there exists $C_T>0$, independent of $\Delta$ and $N$, such that
\begin{equation}\label{eq:mainrate}
\sup_{1\le i\le N}\sup_{0\le t\le T}\E\big|X_t^i-\Zt_t^i\big|\le C_T  \begin{cases}
N^{-1/2}+ \frac{1}{\log(1/\Delta)},&\beta=0\\
 N^{-1/2}+\Delta^\beta,&\beta\in(0,\frac{1}{2}]
 \end{cases}
\end{equation}
Moreover, by Lemma \ref{lem:empirical},
\begin{equation}\label{eq:mainmeasure}
\sup_{0\le t\le T}\E\big[W_1(\mu_t,\hat\mu_t^Z)\big]\le C_T\begin{cases}
 N^{-1/2}+ \frac{1}{\log(1/\Delta)},&\beta=0\\
 N^{-1/2}+\Delta^\beta,&\beta\in(0,\frac{1}{2}]
 \end{cases}
\end{equation}
\end{theorem}

To prove this convergence theorem, we need the following Yamada-Watanabe approximation (see e.g.\ \cite{GR})

For $\delta>1$ and $\varepsilon\in(0,1)$, let
$\psi_{\delta,\varepsilon}:\R\to\R_+$ be continuous with support in
$[\varepsilon/\delta,\varepsilon]$,
$\int_{\varepsilon/\delta}^{\varepsilon}\psi_{\delta,\varepsilon}(z)dz=1$
and $0\le\psi_{\delta,\varepsilon}(z)\le 2/(z\log\delta)$, and set
$\phi_{\delta,\varepsilon}(x)=\int_0^{|x|}\int_0^y
\psi_{\delta,\varepsilon}(z)\,dz\,dy$. Then:
\begin{itemize}
    \item[(P1)] $|x|\le\varepsilon+\phi_{\delta,\varepsilon}(x)$;
    \item[(P2)] $|\phi_{\delta,\varepsilon}'(x)|\le1$ and
    $\phi_{\delta,\varepsilon}'(x)=(x/|x|)\phi_{\delta,\varepsilon}'(|x|)$
    with $0\le\phi_{\delta,\varepsilon}'(|x|)\le1$;
    \item[(P3)]
    $\phi_{\delta,\varepsilon}''(x)=\psi_{\delta,\varepsilon}(|x|)
    \le\frac{2}{|x|\log\delta}\mathbf{1}_{[\varepsilon/\delta,\varepsilon]}(|x|)
    \le\frac{2\delta}{\varepsilon\log\delta}$;
    \item[(P4)] If $a:\R\to\R$ is one-sided Lipschitz with constant
    $\lambda$, i.e.\ $(x-y)(a(x)-a(y))\le\lambda|x-y|^2$, then
    $\phi_{\delta,\varepsilon}'(x-y)(a(x)-a(y))\le\lambda|x-y|$.
\end{itemize}

 Define
\begin{equation*}
e_t^i:=X_t^i-\Zt_t^i,\qquad t\ge0, i=1,\cdots, N.
\end{equation*}

Apply It\^o's formula to $\phi_{\delta,\varepsilon}(e^i_t)$, denote $\phi=\phi_{\delta,\varepsilon}$ for brevity, using \eqref{eq:NIPS} and Lemma \ref{lem:interp}:

\begin{equation}\label{eq:ito}
d\phi(e^i_t)=\phi'(e^i_t)\big(D^i_1(t)+D^i_2(t)\big)\,dt+\tfrac12\phi''(e^i_t)D^i_3(t)\,dt+\phi'(e^i_t)D^i_4(t)\,dW^i_t,
\end{equation}
where
\begin{align}
D^i_1(t)&:=f_1(X^i_t,\mu_t)-\Df_{t-\kappa(t)}(Z^{i,\mathrm{EM}}_t)\,f_1(Z^i_{\kappa(t)},\hat\mu_{\kappa(t)}^Z),\label{eq:D1}\\
D^i_2(t)&:=f_2(X^i_t)-f_2(\Zt^i_t)-\tfrac12\Df^2_{t-\kappa(t)}(Z^{i,\mathrm{EM}}_t)\,g(Z^i_{\kappa(t)})^2,\label{eq:D2}\\
D^i_3(t)&:=\big(g(X^i_t)-\Df_{t-\kappa(t)}(Z^{i,\mathrm{EM}}_t)\,g(Z^i_{\kappa(t)})\big)^2,\label{eq:D3}\\
D^i_4(t)&:=g(X^i_t)-\Df_{t-\kappa(t)}(Z^{i,\mathrm{EM}}_t)\,g(Z^i_{\kappa(t)}).\label{eq:D4}
\end{align}
The term $D^i_1(t)$ is decomposed further as
\begin{equation}\label{eq:D1split}
\aligned D^i_1(t)&=\underbrace{f_1(X^i_t,\mu_t)-f_1(\Zt^i_t,\hat\mu_t^Z)}_{=:D^i_{11}(t)}\\&
\quad+\underbrace{f_1(\Zt^i_t,\hat\mu_t^Z)-\Df_{t-\kappa(t)}(Z^{i,\mathrm{EM}}_t)f_1(Z^i_{\kappa(t)},\hat\mu_{\kappa(t)}^Z)}_{=:D^i_{12}(t)}.\endaligned
\end{equation}
where
\begin{equation}\label{eq:muZt}
\hat\mu_t^Z:=\frac1N\sum_{i=1}^N\delta_{\Zt_t^i},\qquad t\in[0,T],
\end{equation}
the empirical measure of the interpolation.

The following lemma collects the
estimations of the individual terms in \eqref{eq:ito}.

\begin{lemma}\label{lem:drift}
Let Assumptions \ref{ass:A1}-\ref{ass:A3} hold and let
$\phi=\phi_{\delta,\varepsilon}$ with $\delta>1$,
$\varepsilon\in(0,1)$. If $\E |X_0|^{2(l+1)}<\infty,$ then there exists $C_T>0$, independent of
$\Delta$, $N$, $\delta$ and $\varepsilon$, such that for every
$t\in[0,T]$,
\begin{align}
\E\bigl[\phi'(e^i_t)D^i_{11}(t)\bigr]
&\le L\,\E\phi(e^i_t)+L\,\E\Wone(\mu_t,\hat\mu_t^Z)+L\varepsilon,
\label{eq:eD11}\\
\E\bigl[\phi'(e^i_t)D^i_{12}(t)\bigr]
&\le C_T\Delta^{1/2},
\label{eq:eD12}\\
\E\bigl[\phi'(e^i_t)\bigl(f_2(X^i_t)-f_2(\Zt^i_t)\bigr)\bigr]
&\le\lambda\,\E\phi(e^i_t)+\lambda\varepsilon,
\label{eq:eD21}\\
\E\bigl[|\phi'(e^i_t)\Df^2_{t-\kappa(t)}(Z^{i,\mathrm{EM}}_t)\bar B^i(t)^2|\bigr]
&\le C_T\Delta,
\label{eq:eD22}\\
\E\bigl[\phi''(e^i_t)D^i_3(t)\bigr]
&\le\frac{2\rho^2\,\varepsilon^{2\beta}}{\log\delta}
+\frac{C_T\,\delta\,\Delta^{1/2+\beta}}{\varepsilon\log\delta}.
\label{eq:eD3}
\end{align}
\end{lemma}

\begin{proof}
We repeatedly use (P1)--(P4), Lemma \ref{lem:flow}, Lemma
\ref{lem:moments} and Lemma \ref{lem:increments}; $C_T$ may change from
line to line. Set $s=t-\kappa(t)\le\Delta$ and $\kappa=\kappa(t)$.

\emph{Estimate \eqref{eq:eD11}.} By (P2) and Assumption \ref{ass:A1},
\begin{align*}
|\phi'(e^i_t)D^i_{11}(t)|
&\le |f_1(X^i_t,\mu_t)-f_1(\Zt^i_t,\hat\mu_t^Z|\\
&\le  L|e^i_t|+L\Wone(\mu_t,\hat\mu_t^Z),
\end{align*}
and \eqref{eq:eD11} follows from (P1), i.e.\
$|e^i_t|\le\phi(e^i_t)+\varepsilon$.

\emph{Estimate \eqref{eq:eD12}.} Decompose $D^i_{12}=D_{12}^{(a)}+D_{12}^{(b)}$
with
\begin{equation*}
D_{12}^{(a)}:=f_1(\Zt^i_t,\hat\mu_t^Z)-f_1(Z^i_{\kappa},\hat\mu_{\kappa}^Z),
\qquad
D_{12}^{(b)}:=\bigl(1-\Df_{s}(Z^{i,\mathrm{EM}}_t)\bigr)\,
f_1(Z^i_{\kappa},\hat\mu_{\kappa}^Z).
\end{equation*}
By Assumption \ref{ass:A1}, $\Wone(\hat\mu_t^Z,\hat\mu_{\kappa}^Z)\le\frac1N\sum_{i=1}^N|\Zt_t^i-Z_{\kappa}^i|$ and Lemma \ref{lem:increments},
\begin{equation*}
\E|D_{12}^{(a)}|
\le L\,\E|\Zt^i_t-Z^i_{\kappa}|+\frac{L}{N}\sum_{i=1}^N
\E|\Zt_t^i-Z_{\kappa}^i|
\le 2L\sup_{1\le i\le N}\E|\Zt_t^i-Z_{\kappa}^i|
\le C_T\Delta^{1/2}.
\end{equation*}
By (F4) and the
linear growth of $f_1$,
\begin{equation*}
\begin{aligned}
\E|D_{12}^{(b)}|
&\le C\Delta\,\E\bigl[(1+|Z^{i,\mathrm{EM}}_t|^{l})
\bigl|f_1(Z^i_{\kappa},\hat\mu_{\kappa}^Z)\bigr|\bigr]\\
&\le C\Delta\,\sqrt{\E(1+|Z^{i,\mathrm{EM}}_t|^{2l})}
\sqrt{\E|f_1(Z^i_{\kappa},\hat\mu_{\kappa}^Z)|^2}\\
&\le C_T\Delta.
\end{aligned}
\end{equation*}
where we have used Lemma \ref{lem:moments} and
$\E\|\hat\mu_{\kappa}^Z\|_1^2\le \E|Z_{\kappa}^i|^2$. Since
$|\phi'|\le1$, \eqref{eq:eD12} follows.

\emph{Estimate \eqref{eq:eD21}.} This is (P4) applied to $f_2$
(Assumption \ref{ass:A2}), together with (P1).

\emph{Estimate \eqref{eq:eD22}.} By (F4),
$|\Df^2_{s}(Z^{i,\mathrm{EM}}_t)|\le C\Delta(1+|Z^{i,\mathrm{EM}}_t|^{2l})$, and by Assumption
\ref{ass:A3}, $|\bar B^i(t)|^2=|g(Z^i_{\kappa})|^2\le
2\rho'^2(1+|Z^i_{\kappa}|^2)$ where $\rho'=\rho\vee|g(0)|$. Since $|\phi'|\le1$ and, by Young's
inequality,
\begin{equation*}
(1+|Z^{i,\mathrm{EM}}_t|^{2l})(1+|Z^i_{\kappa}|^2)
\le C\bigl(1+|Z^{i,\mathrm{EM}}_t|^{2l+2}+|Z^i_{\kappa}|^{2l+2}\bigr),
\end{equation*}
Lemma \ref{lem:moments} yields
\begin{equation*}
\E\bigl|\phi'(e^i_t)\Df^2_s(Z^{i,\mathrm{EM}}_t)\bar B^i(t)^2\bigr|
\le C\Delta\,\E\bigl[(1+|Z^{i,\mathrm{EM}}_t|^{2l})(1+|Z^i_{\kappa}|^2)\bigr]
\le C_T\Delta.
\end{equation*}

\emph{Estimate \eqref{eq:eD3}.} Write
$D^i_3(t)\le2D_{31}(t)+2D_{32}(t)$ with
\begin{equation*}
D_{31}(t):=\bigl(g(X^i_t)-g(\Zt^i_t)\bigr)^2,\qquad
D_{32}(t):=\bigl(g(\Zt^i_t)-\Df_s(Z^{i,\mathrm{EM}}_t)\,g(Z^i_{\kappa})\bigr)^2.
\end{equation*}

For $D_{31}$, Assumption \ref{ass:A3} and (P3) imply
\begin{equation}\label{eq:D31est}
\phi''(e^i_t)D_{31}(t)
\le\rho^2\,\phi''(e^i_t)\,|e^i_t|^{1+2\beta}
\le\frac{2\rho^2\,|e^i_t|^{2\beta}}{\log\delta}\,
\mathbf 1_{[\varepsilon/\delta,\varepsilon]}(|e^i_t|)
\le\frac{2\rho^2\,\varepsilon^{2\beta}}{\log\delta}.
\end{equation}

For $D_{32}$, decompose
\begin{equation*}
D_{32}(t)\le2\bigl(g(\Zt^i_t)-g(Z^i_{\kappa})\bigr)^2
+2\bigl(1-\Df_s(Z^{i,\mathrm{EM}}_t)\bigr)^2g(Z^i_{\kappa})^2.
\end{equation*}
By Assumption \ref{ass:A3} and Lemma \ref{lem:increments},
\begin{equation*}
\E|g(\Zt^i_t)-g(Z^i_{\kappa})|^2
\le\rho^2\,\E|\Zt^i_t-Z^i_{\kappa}|^{1+2\beta}
\le C_T\Delta^{1/2+\beta}.
\end{equation*}
(F4) in Lemma \ref{lem:flow}, Young's inequality as above and Lemma \ref{lem:moments} implies
\begin{equation*}
\E\bigl[\bigl(1-\Df_s(Z^{i,\mathrm{EM}}_t)\bigr)^2g(Z^i_{\kappa})^2\bigr]
\le C\Delta^2\,\E\bigl[(1+|Z^{i,\mathrm{EM}}_t|^{2l})(1+|Z^i_{\kappa}|^2)\bigr]
\le C_T\Delta^2.
\end{equation*}
Hence $\E D_{32}(t)\le C_T\Delta^{1/2+\beta}$, and the crude bound
$\phi''\le\frac{2\delta}{\varepsilon\log\delta}$ from (P3) gives
\begin{equation*}
\E\bigl[\phi''(e^i_t)D_{32}(t)\bigr]
\le\frac{2\delta}{\varepsilon\log\delta}\,\E D_{32}(t)
\le\frac{C_T\,\delta\,\Delta^{1/2+\beta}}{\varepsilon\log\delta}.
\end{equation*}
Combining the two pieces yields \eqref{eq:eD3}.
\end{proof}

\begin{proof}[Proof of Theorem \ref{thm:main}]

Integrating \eqref{eq:ito} on $[0,t]$,
taking expectations and Lemma \ref{lem:drift} imply
\begin{equation}\label{eq:closure}
\begin{aligned}
\E\phi(e^i_t)
&\le(L+\lambda)\int_0^t\bigl(\E\phi(e^i_s)+\varepsilon\bigr)ds
+L\int_0^t\E\Wone(\mu_s,\hat\mu_s^Z)\,ds\\
&\quad+C_T\Bigl(\Delta^{1/2}+\frac{\varepsilon^{2\beta}}{\log\delta}
+\frac{\delta\,\Delta^{1/2+\beta}}{\varepsilon\log\delta}\Bigr).
\end{aligned}
\end{equation}

Moreover, by the triangle inequality, Lemma
\ref{lem:empirical}, the coupling bound
$$\Wone(\bar\mu_s^{X,N},\hat\mu_s^Z)\le\frac1N\sum_{i=1}^N|e_s^i|$$ and
(P1), we have
\begin{equation}\label{eq:measplit}
\E\Wone(\mu_s,\hat\mu_s^Z)
\le C_TN^{-1/2}+\frac1N\sum_{i=1}^N\E|e_s^i|
\le C_TN^{-1/2}+\sup_{1\le i\le N}\E\phi(e_s^i)+\varepsilon.
\end{equation}
Combining \eqref{eq:closure} and
\eqref{eq:measplit}, Gronwall's lemma yields for every $t\in[0,T]$,
\begin{equation}\label{eq:measplit1}
 \sup_{1\le i\le N}\E\phi(e_t^i)\le C_T\Bigl(N^{-1/2}+\varepsilon+\Delta^{1/2}
 +\frac{\varepsilon^{2\beta}}{\log\delta}
 +\frac{\delta\Delta^{1/2+\beta}}{\varepsilon\log\delta}\Bigr).
\end{equation}

For $\beta=0$ the estimate \eqref{eq:measplit1} reads
\begin{equation}\label{eq:measplit2}
 \sup_{1\le i\le N}\E\phi(e_t^i)\le C_T\Bigl(N^{-1/2}+\varepsilon+\Delta^{1/2}
 +\frac{1}{\log\delta}
 +\frac{\delta\Delta^{1/2}}{\varepsilon\log\delta}\Bigr).
\end{equation}

 Let
$\varepsilon=\Delta^{1/4}$ and $\delta=\Delta^{-1/4}$ (so
$\log\delta=\frac14\log(1/\Delta)$),
\begin{equation*}
\begin{aligned}
\frac{1}{\log\delta}&=\frac{4}{\log(1/\Delta)},\qquad
\varepsilon=\Delta^{1/4}\le\frac{C}{\log(1/\Delta)},\\
\frac{\delta\,\Delta^{1/2}}{\varepsilon\log\delta}
&=\frac{4}{\log(1/\Delta)},\qquad
\Delta^{1/2}\le\frac{C}{\log(1/\Delta)}.
\end{aligned}
\end{equation*}
and Gronwall's lemma gives, for every $t\in[0,T]$,
\[
  \sup_{1\le i\le N}\sup_{0\le t\le T}\E\phi(e_t^i)\le C_T\bigl(N^{-1/2}+1/\log(1/\Delta)\bigr).
\]
Hence
\begin{equation*}
\sup_{1\le i\le N}\sup_{0\le t\le T}\E|e_t^i|
\le \sup_{1\le i\le N}\sup_{0\le t\le T}\E\phi(e_t^i)+\varepsilon
\le C_T\bigl(N^{-1/2}+1/\log(1/\Delta)\bigr).
\end{equation*}

For $0<\beta\le\frac12$, take $\delta=2$ (so
$\log\delta=\log2$) and $\varepsilon=\Delta^{1/2}$. Moreover, because $0<\Delta<1$, $\varepsilon=\Delta^{1/2}\le\Delta^\beta$. Then
\[
 \sup_{1\le i\le N}\E\phi(e_t^i)\le C_T\Bigl(N^{-1/2}+\Delta^{\beta}\Bigr).
\]

Consequently,
\begin{equation*}
\sup_{0\le t\le T}\sup_{1\le i\le N}\E|e_t^i|
\le \sup_{1\le i\le N}\sup_{0\le t\le T}\E\phi(e_t^i)+\varepsilon
\le C_T\bigl(\Delta^\beta+N^{-1/2}\bigr).
\end{equation*}
which is \eqref{eq:mainrate}. The measure estimate
\eqref{eq:mainmeasure} follows from \eqref{eq:measplit} and Lemma
\ref{lem:empirical}.

\end{proof}

\begin{theorem}\label{thm:pathsup1}
Let Assumptions \ref{ass:A1}-\ref{ass:A3} hold. For $p\ge1$ and $q>(2p\beta)\vee(1+2\beta)$,
if $\E |X_0|^{(p\vee q)(l+1)}<\infty,$ then for $0<\Delta\le e^{-1}$, it follows that

\begin{equation}\label{eq:pathsup-rate1}
\aligned &\quad\sup_{1\le i\le N} \E\left(\sup_{0\le t\le T}
 |X_t^i-\widetilde Z_t^i|^p\right)\\&\le C\left\{
\begin{aligned}
&\left(N^{-1/2}+\frac1{\log(1/\Delta)}\right)^{\frac{p}{2}\wedge1}, & &p\ge1, \beta=0\\
&\left(N^{-1/2}+\Delta^\beta\right)^{(2\beta)\vee \frac{q-1-2\beta}{2(q-1)}},  & &p=1, 0<\beta\le\frac{1}{2} \\
&\left(N^{-1/2}+\Delta^\beta\right)^{\frac{p(q-1-2\beta)}{2(q-1)}}, & &1<p\le2, 0<\beta\le\frac{1}{2} \\
&\left(N^{-1/2}+\Delta^\beta\right)^{2p\beta}, & &p\ge2, 0< \beta\le \frac{1}{2p} \\
&\left(N^{-1/2}+\Delta^\beta\right)^\frac{q-2p\beta}{q-1}, & &p\ge2, \frac{1}{2p}<\beta\le\frac{1}{2}
\end{aligned}
\right.\endaligned
\end{equation}
where $C>0$ is a constant which is independent of $\Delta$ and $N.$

Specially,
\begin{equation}\label{eq:pathsup-rate}
\sup_{1\le i\le N} \mathbb E\!\left[\sup_{0\le t\le T}
 |X_t^i-\widetilde Z_t^i|\right]
 \le C\begin{cases}
  \bigl(N^{-1/2}+\frac1{\log(1/\Delta)}\bigr)^{1/2},
     &\beta=0\\
 \bigl(N^{-1/2}+\Delta^\beta\bigr)^{(2\beta)\vee\frac{q-1-2\beta}{2(q-1)}},
     &0<\beta\le\frac12
 \end{cases}
\end{equation}
\begin{equation}\label{eq:pathsup-rate'}
\sup_{1\le i\le N} \mathbb E\!\left[\sup_{0\le t\le T}
 |X_t^i-\widetilde Z_t^i|^2\right]
 \le C\begin{cases}
  \bigl(N^{-1/2}+\frac1{\log(1/\Delta)}\bigr),
     &\beta=0\\
 \bigl(N^{-1/2}+\Delta^\beta\bigr)^{(4\beta)\vee\frac{q-1-2\beta}{q-1}},
     &0<\beta\le\frac14,\\
 \left(N^{-1/2}+\Delta^\beta\right)^\frac{q-4\beta}{q-1}, & \frac{1}{4}<\beta\le\frac{1}{2}

 \end{cases}
\end{equation}

\end{theorem}

\begin{proof}
By It\^o's formula and Lemma \ref{lem:drift}, taking the $p$-th power and expection on both sides ($p\ge1$), it follows that
\begin{equation}\label{eq:closure1'}
\begin{aligned}
&\quad\E\sup_{0\le s\le t}\phi(e^i_s)^p\\
&\le C\E\int_0^t\left(|e^i_s|^p
+\Wone^p(\mu_s,\hat\mu_s^Z)\right)\,ds\\
& \quad+ C\E\int_0^t\left(|\tilde{Z}_s^i-Z_{\kappa(s)}^i|^p
+\frac{1}{N}\sum_{i=1}^N|\tilde{Z}_s^i-Z_{\kappa(s)}^i|^p\right)\,ds\\&
\quad+C\E\int_0^t\left(\phi(e^i_s)^p+\varepsilon^p\right)\,ds\\&
\quad+C\Delta^p\E\int_0^t\left(1+|Z^{i,\mathrm{EM}}_t|^{p(l+1)}+|Z^i_{\kappa}|^{p(l+1)}\right)\,ds
\\&\quad+C\frac{\varepsilon^{2p\beta}}{\log^p\delta}+C\frac{\delta^p}{\varepsilon^p\log^p\delta}\E\left[\int_0^t\left(|\tilde{Z}_s^i-Z_{\kappa(s)}^i|^{p(1+2\beta)}
\right)\,ds\right.\\&
\left.\quad+C\Delta^{2p}\E\int_0^t\left(1+|Z^{i,\mathrm{EM}}_t|^{2p(l+1)}+|Z^i_{\kappa}|^{2p(l+1)}\right)\,ds\right]+C\E\sup_{0\le s\le t}M^p_s
\end{aligned}
\end{equation}
where $M_t=\int_0^t\phi'(e^i_s)D^i_4(s)\,dW^i_s.$

Since
$$\aligned \E\Wone^p(\mu_s,\hat\mu_s^Z)&\le \E(\Wone(\mu_s,\bar{\mu}_s^{X,N})+\Wone(\bar{\mu}_s^{X,N},\hat\mu_s^Z))^p\\&
\le C_p(N^{-\frac{p}{2}}+\E|e_s^i|^p)\endaligned$$
and the Burkholder--Davis--Gundy inequality implies
\begin{equation}\label{eq:closure2'}\aligned \E\sup_{0\le s\le t}M^p_s&\le C\E \left(\int_0^t|\phi'(e^i_s)|^2D^i_3(s)ds\right)^\frac{p}{2}\\&
\le C\E \left(\int_0^t\left(2\rho^2|e^i_s|^{1+2\beta}+4\rho^2|\Zt^i_s-Z^i_{\kappa(s)}|^{1+2\beta}\right.\right.\\&
\qquad+\left.\left.\Delta^2\bigl[(1+|Z^{i,\mathrm{EM}}_s|^{2l})(1+|Z^i_{\kappa(s)}|^2)\bigr]\right)ds\right)^\frac{p}{2}\\&
\le C\left(\E\left(\int_0^t|e^i_s|^{1+2\beta}ds\right)^\frac{p}{2}+\E\left(\int_0^t|\Zt^i_s-Z^i_{\kappa(s)}|^{1+2\beta}ds\right)^\frac{p}{2}\right.\\&
\qquad+\left.\Delta^p\E\left(\int_0^t\bigl[(1+|Z^{i,\mathrm{EM}}_s|^{2(l+1)}+|Z^i_{\kappa(s)}|^{2(l+1)}\bigr]ds\right)^\frac{p}{2}\right)\\&
\le C \left(\E\left(\int_0^t|e^i_s|^{1+2\beta}ds\right)^\frac{p}{2}+\Delta^{\frac{p(1+2\beta)}{4}}\right).
\endaligned
\end{equation}

Then
$$\aligned \E\sup_{0\le s\le t}\phi(e^i_s)^p&\le C\int_0^t\E\phi(e^i_s)^p\,ds+C\left(\varepsilon^p+N^{-\frac{p}{2}}+\Delta^\frac{p}{2}+\frac{\varepsilon^{2p\beta}}{\log^p\delta}+\Delta^p\right.\\&
\left.\quad+\frac{\delta^p}{\varepsilon^p\log^p\delta}\Delta^\frac{p(1+2\beta)}{2}\right)+C \left(\E\left(\int_0^t|e^i_s|^{1+2\beta}ds\right)^\frac{p}{2}+\Delta^{\frac{p(1+2\beta)}{4}}\right).\endaligned
$$

Gronwall's lemma implies
\begin{equation}\label{gronwall}\aligned \E\sup_{0\le s\le t}\phi(e^i_s)^p&\le C\left(\varepsilon^p+N^{-\frac{p}{2}}+\Delta^{\frac{p(1+2\beta)}{4}}+\frac{\varepsilon^{2p\beta}}{\log^p\delta}+\frac{\delta^p}{\varepsilon^p\log^p\delta}\Delta^\frac{p(1+2\beta)}{2}\right)\\&
\quad+C \E\left(\int_0^t|e^i_s|^{1+2\beta}ds\right)^\frac{p}{2}.\endaligned \end{equation}

If $\beta=0, 1\le p\le 2,$ then we have
$$\aligned \E\sup_{0\le s\le t}\phi(e^i_s)^p&\le C\left(\varepsilon^p+N^{-\frac{p}{2}}+\Delta^{\frac{p}{4}}+\frac{1}{\log^p\delta}+\frac{\delta^p}{\varepsilon^p\log^p\delta}\Delta^\frac{p}{2}\right)\\&
\quad+C\left(\int_0^t\E|e^i_s|ds\right)^\frac{p}{2}\\&
\le C\left(\varepsilon^p+N^{-\frac{p}{2}}+\Delta^\frac{p}{4}+\frac{1}{\log^p\delta}+\frac{\delta^p}{\varepsilon^p\log^p\delta}\Delta^\frac{p}{2}\right) \\&
\quad +C\left(N^{-1/2}+\frac1{\log(1/\Delta)}\right)^\frac{p}{2}.\endaligned $$

If $\beta=0,  p>2,$ then
$$\aligned \E\sup_{0\le s\le t}\phi(e^i_s)^p&\le C\left(\varepsilon^p+N^{-\frac{p}{2}}+\Delta^{\frac{p}{4}}+\frac{1}{\log^p\delta}+\frac{\delta^p}{\varepsilon^p\log^p\delta}\Delta^\frac{p}{2}\right)\\&
\quad+C\left(\int_0^t\E|e^i_s|^\frac{p}{2}ds\right)\\&
\le C\left(\varepsilon^p+N^{-\frac{p}{2}}+\Delta^\frac{p}{4}+\frac{1}{\log^p\delta}+\frac{\delta^p}{\varepsilon^p\log^p\delta}\Delta^\frac{p}{2}\right) \\&
\quad +C\left(\int_0^t\E|e^i_s|ds+\int_0^t\E|e^i_s|^pds\right)\\&
\le C\left(\varepsilon^p+N^{-\frac{p}{2}}+\Delta^\frac{p}{4}+\frac{1}{\log^p\delta}+\frac{\delta^p}{\varepsilon^p\log^p\delta}\Delta^\frac{p}{2}\right) \\&
\quad +C\left(N^{-1/2}+\frac1{\log(1/\Delta)}\right)+C\int_0^t\E\phi(e^i_s)^pds.\endaligned $$

By using Gronwall's lemma again, it follows that
$$\aligned \E\sup_{0\le s\le t}\phi(e^i_s)^p&
\le C\left(\varepsilon^p+N^{-\frac{p}{2}}+\Delta^\frac{p}{4}+\frac{1}{\log^p\delta}+\frac{\delta^p}{\varepsilon^p\log^p\delta}\Delta^\frac{p}{2}\right) \\&
\quad +C\left(N^{-1/2}+\frac1{\log(1/\Delta)}\right).\endaligned $$

Therefore, if $\beta=0,$ then for any $p\ge1,$ we have
$$\aligned \E\sup_{0\le s\le t}|e^i_s|^p&\le C\E\sup_{0\le s\le t}\phi(e^i_s)^p+C\varepsilon^p\\&
\le  C\left(\varepsilon^p+N^{-\frac{p}{2}}+\Delta^{\frac{p}{2}}
+\frac{1}{\log^p\delta}+\frac{\delta^p}{\varepsilon^p\log^p\delta}\Delta^{\frac{p}{2}}\right)\\&
\quad+C\left(N^{-1/2}+\frac1{\log(1/\Delta)}\right)^{\frac{p}{2}\wedge1}.\endaligned$$

Letting $\varepsilon=\Delta^{1/4}$ and $\delta=\Delta^{-1/4}$, we have
\begin{equation}\label{log}\aligned \E\sup_{0\le s\le t}|e^i_s|^p&
\le C\left(N^{-1/2}+\frac1{\log(1/\Delta)}\right)^{\frac{p}{2}\wedge1}.\endaligned\end{equation}

Now let us consider deal with $0<\beta\le\frac{1}{2}$. The value of $p$ can be divided into the following four cases:

\textbf{Case 1}: $p=1$. Then for any $0<\beta\le\frac{1}{2}$,
$$\aligned \E\left(\int_0^t|e^i_s|^{1+2\beta}ds\right)^\frac{1}{2}&\le \left(\int_0^t\E|e^i_s|^{1+2\beta}ds\right)^\frac{1}{2}\\&
\le \left(\int_0^t(\E|e^i_s|)^\frac{q-1-2\beta}{q-1}(\E|e_t^i|^q)^\frac{2\beta}{q-1}ds\right)^\frac{1}{2}\\&
\le C \left(N^{-1/2}+\Delta^\beta\right)^{\frac{q-1-2\beta}{2(q-1)}}.\endaligned$$

On the other hand,
$$\aligned \E\left(\int_0^t|e^i_s|^{1+2\beta}ds\right)^\frac{1}{2}&\le \E\left(\sup_{0\le s\le t}|e^i_s|^\frac{1}{2}\left(\int_0^t|e^i_s|^{2\beta}ds\right)^\frac{1}{2}\right)\\&
\le \varepsilon'\E\sup_{0\le s\le t}|e^i_s|+\frac{1}{4\varepsilon'}\left(\int_0^t\E|e^i_s|^{2\beta}ds\right)\\&
\le \varepsilon'\E\sup_{0\le s\le t}|e^i_s|+C\left(N^{-1/2}+\Delta^\beta\right)^{2\beta}.\endaligned$$

Thus
$$\aligned \E\left(\int_0^t|e^i_s|^{1+2\beta}ds\right)^\frac{1}{2}&\le \varepsilon'\E\sup_{0\le s\le t}|e^i_s|+C\left(N^{-1/2}+\Delta^\beta\right)^{(2\beta)\vee \frac{q-1-2\beta}{2(q-1)}}.\endaligned$$

\textbf{Case 2}: $1< p\le2$. Similar to Case 1, for any $0<\beta\le\frac{1}{2}$, it follows that
$$\aligned \E\left(\int_0^t|e^i_s|^{1+2\beta}ds\right)^\frac{p}{2}&
\le C \left(N^{-1/2}+\Delta^\beta\right)^{\frac{p(q-1-2\beta)}{2(q-1)}}.\endaligned$$

\textbf{Case 3}: $p\ge2$ and $\beta\le \frac{1}{2p}$. We have
$$\aligned \E\left(\int_0^t|e^i_s|^{1+2\beta}ds\right)^\frac{p}{2}&\le \E\left(\sup_{0\le s\le t}|e^i_s|^\frac{p}{2}\left(\int_0^t|e^i_s|^{2\beta}ds\right)^\frac{p}{2}\right)\\&
\le \varepsilon''\E\sup_{0\le s\le t}|e^i_s|^p+\frac{1}{4\varepsilon''}\left(\int_0^t\E|e^i_s|^{2p\beta}ds\right)\\&
\le \varepsilon''\E\sup_{0\le s\le t}|e^i_s|^p+C \left(N^{-1/2}+\Delta^\beta\right)^{2p\beta}.\endaligned$$

\textbf{Case 4}: $p\ge2$ and $ \frac{1}{2p}<\beta\le\frac{1}{2}$. In this case, $1<2p\beta\le p$. Then
$$\aligned \E\left(\int_0^t|e^i_s|^{1+2\beta}ds\right)^\frac{p}{2}&\le \varepsilon'''\E\sup_{0\le s\le t}|e^i_s|^p+\frac{1}{4\varepsilon'''}\left(\int_0^t\E|e^i_s|^{2p\beta}ds\right)\\&
\le\varepsilon'''\E\sup_{0\le s\le t}|e^i_s|^p+C\int_0^t(\E|e^i_s|)^\frac{q-2p\beta}{q-1}(\E|e_t^i|^q)^\frac{2p\beta-1}{q-1}ds\\&
\le \varepsilon'''\E\sup_{0\le s\le t}|e^i_s|^p+C \left(N^{-1/2}+\Delta^\beta\right)^\frac{q-2p\beta}{q-1}.\endaligned$$

Choosing $\varepsilon', \varepsilon''$ and $\varepsilon'''$ sufficiently small, then \eqref{gronwall} yields
$$\aligned \E\sup_{0\le s\le t}|e^i_s|^p&\le\E\sup_{0\le s\le t}\phi(e^i_s)^p+\varepsilon^p\\& \le C\left(\varepsilon^p+N^{-\frac{p}{2}}+\Delta^{\frac{p(1+2\beta)}{4}}+\frac{\varepsilon^{2p\beta}}{\log^p\delta}+\frac{\delta^p}{\varepsilon^p\log^p\delta}\Delta^{\frac{p(1+2\beta)}{2}}\right)\\&
\quad
+C\left\{
\begin{aligned}
&\left(N^{-1/2}+\Delta^\beta\right)^{(2\beta)\vee \frac{q-1-2\beta}{2(q-1)}},  &p=1, 0<\beta\le\frac{1}{2} \\
&\left(N^{-1/2}+\Delta^\beta\right)^{\frac{p(q-1-2\beta)}{2(q-1)}},  &1\le p\le2, 0<\beta\le\frac{1}{2} \\
&\left(N^{-1/2}+\Delta^\beta\right)^{2p\beta}, & p\ge2, 0< \beta\le \frac{1}{2p} \\
&\left(N^{-1/2}+\Delta^\beta\right)^\frac{q-2p\beta}{q-1}, & p\ge2, \frac{1}{2p}<\beta\le\frac{1}{2}
\end{aligned}
\right.\endaligned$$

By letting $\delta=2$, $\varepsilon=\Delta^{1/2}$ and eliminating the higher order term, it follows that
\begin{equation}\label{yizhi}
\aligned \E\sup_{0\le s\le t}|e^i_s|^p&\le C\left\{
\begin{aligned}
&\left(N^{-1/2}+\Delta^\beta\right)^{(2\beta)\vee \frac{q-1-2\beta}{2(q-1)}},  &p=1, 0<\beta\le\frac{1}{2} \\
&\left(N^{-1/2}+\Delta^\beta\right)^{\frac{p(q-1-2\beta)}{2(q-1)}},  &1\le p\le2, 0<\beta\le\frac{1}{2} \\
&\left(N^{-1/2}+\Delta^\beta\right)^{2p\beta}, & p\ge2, 0< \beta\le \frac{1}{2p} \\
&\left(N^{-1/2}+\Delta^\beta\right)^\frac{q-2p\beta}{q-1}, & p\ge2, \frac{1}{2p}<\beta\le\frac{1}{2}
\end{aligned}
\right.\endaligned\end{equation}

So \eqref{eq:pathsup-rate1} is a direct consequence of \eqref{log} and \eqref{yizhi}.

We then complete the proof.
\end{proof}

\begin{remark}\label{zhu}
In \cite{LiuShiWu}, Corollary 3.3, authors obtain the $L^1$ and $L^2$ uniform convergence rates of tamed EM method. Comparing with their results, under slightly stronger condition Assumption \ref{ass:A2}, the $L^1$ convergence rate of our splitting Euler method strictly better than that of (3.9) since (see \eqref{eq:pathsup-rate}). Indeed, for $0<\beta<\frac{1}{4}$, if we choose $q>1+\frac{2\beta}{1-4\beta}$, then $\beta\left[(2\beta)\vee\frac{q-1-2\beta}{2(q-1)}\right]>2\beta^2$. Moreover, $L^2$ convergence rate in $N$ is higher than that of (3.10). for $\beta=0$, or $\beta\in(0,\frac{1}{8}], q>1+4\beta$, or $\beta\in(\frac{1}{8},\frac{1}{4}]$, or $\beta\in(\frac{1}{4},\frac{1}{2}], q>8\beta-1$.
\end{remark}

\section{Stability of exact and numerical solutions}\label{sec:stab}

In this section, we will consider the $L^1$ and $L^2$ exponential stability of exact and the corresponding numerical solution respectively. 
Stability of exact and modified tamed Euler-Maruyama scheme for
McKean--Vlasov equations with H\"older diffusion have recently been studied
in \cite{LiuGaoGuo}. Inspired by this work, we will modify our splitting Euler scheme with exponentially tamed diffusion. The new scheme is defined as follows.

\begin{definition}\label{def:scheme1}
Let $\Delta\in(0,1)$, $t_k=k\Delta$, and
$\hat\mu_k^Z=\frac1N\sum_{j=1}^N\delta_{Z_k^j}$. Set
$Z_0^i=X_0^i$ and, for $k\ge0$,
\begin{equation}\label{eq:scheme1}
    Z_{k+1}^i
    =\flow_\Delta\Bigl(
    Z_k^i+f_1\bigl(Z_k^i,\hat\mu_k^Z\bigr)\Delta
    +g_\Delta\bigl(t_k,Z_k^i\bigr)\,\Delta W_k^i
    \Bigr),
    \qquad
    \Delta W_k^i=W_{t_{k+1}}^i-W_{t_k}^i,
\end{equation}
where $g_\Delta$ is the exponentially tamed diffusion defined
with a fixed parameter $\lambda>0$,
\begin{equation}\label{eq:taming1}
    g_\Delta(t_k,x)
    =\frac{g(x)}{1+\sqrt\Delta\,e^{2\lambda t_k}(1+|g(x)|)}.
\end{equation}
\end{definition}
The Euler-stage process and flow-corrected interpolation are defined the same as \eqref{eq:zem} and \ref{eq:interp} to obtain the new splitting Euler scheme with exponentially tamed diffusion.

We first establish an $L^1$ exponential stability estimate for the exact
solution under strict dissipativity, using the Yamada--Watanabe
approximation from Section \ref{sec:model}. Then we derive an exponential
decay estimate for the second moment of the new numerical scheme \eqref{eq:scheme1} when the diffusion term $g$ is replaced by the
exponentially tamed diffusion \eqref{eq:taming1}.
Throughout this section, the equilibrium and strict dissipativity
condition is
\begin{assumption}\label{ass:S}
$f_1(0,\delta_0)=0$, $f_2(0)=0$, $g(0)=0$, and the dissipativity in
Assumption \ref{ass:A2} holds with vanishing constant and dominates the
Lipschitz constant of $f_1$.
\begin{equation*}
xf_2(x)\le-\lambda'|x|^2.
\end{equation*}
\end{assumption}

Under Assumption \ref{ass:S}, $X_t\equiv0$ is an
equilibrium solution of \eqref{eq:MVSDE}. The differential inequality
used in (F6) follows directly from the strict dissipativity condition;
no claim of uniqueness among all invariant probability measures is
needed for the stability theorem.

\begin{theorem}\label{thm:exact-L1-stab}
Suppose Assumptions \ref{ass:A1}-\ref{ass:A3} and \ref{ass:S} hold with $\lambda'>2L$.
Let $X_t$ be the strong solution of \eqref{eq:MVSDE}. Then, for every
$t\ge0$ and $\beta\in[0,\frac12]$,
\begin{equation}\label{eq:exact-L1-stab}
  \mathbb E|X_t|\le
  e^{-(\lambda'-2L)t}\mathbb E|X_0|.
\end{equation}
\end{theorem}

\begin{proof}
Write $m(t)=\mathbb E|X_t|$ and
$\gamma=\lambda'-2L>0$. Let
$\phi_{\delta,\varepsilon}$ be the Yamada--Watanabe approximation
introduced in Section \ref{sec:model}. Its construction gives
$0\le\phi_{\delta,\varepsilon}(x)\le|x|\le\varepsilon+\phi_{\delta,\varepsilon}(x)$, and
\[
  x\phi'_{\delta,\varepsilon}(x)=|x|\phi'_{\delta,\varepsilon}(|x|)\ge \phi_{\delta,\varepsilon}(|x|)\ge |x|-\varepsilon,
  \quad
  0\le\phi''_{\delta,\varepsilon}(x)
  \le\frac{2}{|x|\log\delta}
  \mathbf 1_{[\varepsilon/\delta,\varepsilon]}(|x|).
\]
The equilibrium condition and Assumption \ref{ass:A1} imply
\[
 |f_1(X_t,\mu_t)|\le L\bigl(|X_t|+W_1(\mu_t,\delta_0)\bigr),
 \qquad \E W_1(\mu_t,\delta_0)=m(t),
\]
so that
$\mathbb E|f_1(X_t,\mu_t)|\le2Lm(t)$. Moreover,
Assumption \ref{ass:S} yields
\[
 \aligned \phi'_{\delta,\varepsilon}(x)f_2(x)&= \frac{\phi'_{\delta,\varepsilon}(x)}{x}xf_2(x)\le -\lambda'|x|^2\frac{\phi'_{\delta,\varepsilon}(x)}{x}=-\lambda'|x|\phi'_{\delta,\varepsilon}(|x|)\\&
 \le-\lambda'(|x|-\varepsilon).\endaligned
\]
Since $g(0)=0$ and Assumption \ref{ass:A3} holds,
$|g(x)|^2\le\rho^2|x|^{1+2\beta}$. Consequently,
\[
\begin{aligned}
\frac12\phi''_{\delta,\varepsilon}(x)|g(x)|^2
&\le \frac{\rho^2}{\log\delta}|x|^{2\beta}
  \mathbf 1_{[\varepsilon/\delta,\varepsilon]}(|x|)\le \frac{\rho^2\varepsilon^{2\beta}}{\log\delta}.
\end{aligned}
\]

Then It\^o's formula implies that for every $t\ge0$,
\[
\begin{aligned}
y(t)=y(0)+\int_0^t\mathbb E\Big[&
\phi'_{\delta,\varepsilon}(X_s)
  \big(f_1(X_s,\mu_s)+f_2(X_s)\big)\\
&+\frac12\phi''_{\delta,\varepsilon}(X_s)|g(X_s)|^2
\Big] ds\\&
\le y(0)+\int_0^t(-\gamma m(s)+c_\varepsilon)ds
\end{aligned}
\]
where $y(t):=\mathbb E\phi_{\delta,\varepsilon}(X_t)$ and
\[
 c_\varepsilon:=\lambda'\varepsilon+
 \frac{\rho^2\varepsilon^{2\beta}}{\log\delta}.
\]
In particular, $y$ is absolutely continuous on every finite interval,
and the integrand estimate gives, for almost every $t\ge0$,
\[
 y'(t)\le -\gamma m(t)+c_\varepsilon
 \le -\gamma y(t)+c_\varepsilon,
\]
where the last inequality uses $y(t)\le m(t)$ and $\gamma>0$.
The differential comparison principle therefore gives
\[
 y(t)\le e^{-\gamma t}y(0)
 +\frac{c_\varepsilon}{\gamma}(1-e^{-\gamma t}).
\]
Using $y(0)\le\mathbb E|X_0|$ and $m(t)\le y(t)+\varepsilon$, we obtain
\[
 m(t)\le e^{-\gamma t}\mathbb E|X_0|+\varepsilon
 +\frac{\lambda'\varepsilon+
       \rho^2\varepsilon^{2\beta}/\log\delta}{\gamma}.
\]
Choose $\delta=\exp(1/\varepsilon)$. Then
$\varepsilon^{2\beta}/\log\delta=\varepsilon^{1+2\beta}\to0$
as $\varepsilon\downarrow0$, including when $\beta=0$.
Letting $\varepsilon\downarrow0$ proves
\eqref{eq:exact-L1-stab}.
\end{proof}

\begin{remark}\label{zhu1}
In \cite{LiuGaoGuo}, Theorem 3.1, the authors obtain the $L^q$ ($0<q\le1$) exponential stability of the MVSDEs under similar conditions. However, the decay rate $\gamma_1$ in \cite{LiuGaoGuo} is smaller than $\rho_1-3\rho_2$. Note that their $\rho_1$ is our $\lambda'$, $\rho_2$ is our $L$. So our decay rate $\lambda'-2L$ is bigger than that of \cite{LiuGaoGuo}. Moreover, our decay does not depend on $\beta.$
\end{remark}

Now let us investigate mean square exponential stability of the splitting Euler scheme with exponentially tamed diffusion.

The pointwise taming bound
\begin{equation}\label{eq:tamingbound}
|g_\Delta(t_k,x)|
=\frac{|g(x)|}{1+\sqrt\Delta\,e^{2\lambda t_k}(1+|g(x)|)}
\le\Delta^{-1/2}e^{-2\lambda t_k}
\end{equation}
is the key to the following result.

\begin{theorem}\label{thm:stab}
Let Assumptions \ref{ass:A1}--\ref{ass:A3} and \ref{ass:S} hold, and
let $Z_k^i$ be the scheme \eqref{eq:scheme1} with the exponentially
tamed diffusion \eqref{eq:taming1}, where $\lambda>0$. Put
$\gamma:=\lambda'-2L>0$ and
\[
  \eta_*:=\min\{2\gamma,4\lambda\}.
\]
Then, for every fixed $\Delta\in(0,1)$ and every
$\eta\in(0,\eta_*)$,
\begin{equation}\label{eq:stab-limit}
  \limsup_{k\to\infty}\frac{\log \sup_{1\le i\le N}\mathbb E|Z_k^i|^2}{t_k}\le-\eta_*.
\end{equation}
That is, for each fixed step size, the numerical second moment
decays exponentially to zero. Moreover,
\begin{equation}\label{eq:stab-limit1}
  \limsup_{k\to\infty}\frac{\log|Z_k^i|}{t_k}\le-\frac{\eta_*}{2}, \quad \textrm{a.s.}
\end{equation}
\end{theorem}

\begin{proof}
Set $a_k:=\sup_{1\le i\le N}\mathbb E|Z_k^i|^2$ and
$Y_k^i:=Z_k^i+f_1(Z_k^i,\hat\mu_k^Z)\Delta
+g_\Delta(t_k,Z_k^i)\Delta W_k^i$. By strict dissipativity Assumption \ref{ass:S}, it is obvious that
\[
 \frac{d}{dt}|y_x(t)|^2
 =2y_x(t)f_2(y_x(t))\le-2\lambda'|y_x(t)|^2.
\]
Then we have $|\flow_\Delta(x)|^2=y_x(\Delta)\le e^{-2\lambda'\Delta}|x|^2$, and hence
\[
\mathbb E|Z_{k+1}^i|^2\le e^{-2\lambda'\Delta}\mathbb E|Y_k^i|^2.
\]
Expanding the square and taking conditional expectation eliminates the
Brownian cross term. Assumption \ref{ass:A1}, the equilibrium condition,
and \[\mathbb E( N^{-1}\sum_{j=1}^N|Z_k^j|)^2\le a_k\] yield
\[
\mathbb E|Z_k^i+\Delta f_1(Z_k^i,\hat\mu_k^Z)|^2
\le(1+4L\Delta+4L^2\Delta^2)a_k.
\]
The taming denominator gives the pointwise estimate
\[
\Delta|g_\Delta(t_k,x)|^2\le e^{-4\lambda t_k}.
\]
Consequently,
\begin{equation}\label{eq:recursion}
 a_{k+1}\le q(\Delta)a_k+e^{-4\lambda t_k}
\end{equation}
where $ q(\Delta)=e^{-2\lambda'\Delta}(1+2L\Delta)^2.$
Since $\log(1+u)\le u$ for $u>-1$,
$q(\Delta)\le e^{-2\gamma\Delta}$. Iterating \eqref{eq:recursion}, it follows that
\[
 a_k\le e^{-2\gamma t_k}a_0+
 \sum_{j=0}^{k-1}e^{-2\gamma(k-1-j)\Delta}e^{-4\lambda j\Delta}.
\]
If $\gamma\ne2\lambda$, the sum equals
\[
 \frac{e^{-2\gamma k\Delta}-e^{-4\lambda k\Delta}}
 {e^{-2\gamma\Delta}-e^{-4\lambda\Delta}}.
\]
For $0<\Delta<1$, the absolute value of the denominator is bounded
below by $c\Delta$, where $c>0$ depends only on the fixed rates.
The numerator is bounded by $e^{-\eta_*t_k}$. Therefore,
\begin{equation}\label{eq:stab}
  \mathbb E|Z_k^i|^2
  \le e^{-2\gamma t_k}\mathbb E|X_0|^2
  +\frac{C}{\Delta}e^{-\eta_*t_k}.
\end{equation}

 If $\gamma=2\lambda$, the sum is
$k e^{-2\gamma(k-1)\Delta}$. So
\begin{equation}\label{eq:stab1}
  \mathbb E|Z_k^i|^2
  \le e^{-2\gamma t_k}\mathbb E|X_0|^2
  +k e^{-2\gamma(k-1)\Delta}.
\end{equation}
Either of \eqref{eq:stab} or \eqref{eq:stab1} means \eqref{eq:stab-limit}.

\eqref{eq:stab-limit1} can be obtained by using \eqref{eq:stab-limit}, Chebyshev inequality and Borel--Cantelli lemma directly.

\end{proof}

\section{Numerical Experiments}
\label{sec:numerical}

In this section, we present numerical experiments to illustrate
the convergence and long-time behavior of the proposed splitting
scheme. We investigate the terminal-time $L^1$ error, the
mean-square decay, and representative numerical sample paths.

\subsection{$L^1$ convergence errors with respect to stepsize}

We first consider the following McKean--Vlasov SDE:
\begin{equation}
\label{eq:convergence-test}
\mathrm{d}X_t
=
\left(0.4\,\mathbb{E}[X_t]-X_t^3\right)\mathrm{d}t
+0.5|X_t|^\alpha\,\mathrm{d}W_t,
\qquad X_0=1.2\xi,
\end{equation}
where $\xi$ is a standard normally distributed random variable.
We consider three H\"older exponents,
$\alpha=0.5,
\alpha=0.75
$ and $\alpha=0.95$, corresponding to $\beta=0$, $\beta=0.25$ and $\beta=0.45$, respectively.

The drift is decomposed as
\[
f_1(x,\mu)
=0.4\int_{\mathbb{R}}y\,\mu(\mathrm{d}y),
\qquad
f_2(x)=-x^3.
\]
The exact flow of the dissipative component is
\[
\Phi_\Delta(y)
=
\frac{y}{\sqrt{1+2\Delta y^2}}.
\]
Consequently, the splitting scheme is given by
\begin{align}
Y_k^i
&=Z_k^i+0.4\overline{Z}_k\Delta
+0.5|Z_k^i|^\alpha\Delta W_k^i,
\quad Z_{k+1}^i
=\Phi_\Delta(Y_k^i),
\end{align}
\eqref{eq:zem} and \eqref{eq:interp}. Where $\overline{Z}_k
=\frac{1}{N}\sum_{j=1}^N Z_k^j.$

We set the terminal time to $T=1$ and use $N=1000$
particles. The reference exact solution is computed with the step size
$\Delta_{\mathrm{ref}}=2^{-15}$, while the coarser step sizes are
\[
\Delta=2^{-6},\,2^{-7},\,2^{-8},\,2^{-9},\,2^{-10},\,2^{-11}.
\]
For each repetition, the Brownian increments on a coarse grid
are obtained by summing the corresponding increments on the
reference grid. The same initial samples and Brownian paths
are used for the two choices of $\alpha$.

We consider two different $L^1$ errors. 
\[
e_1(\Delta)
=
\sup_{0\leq t\leq 1}
\mathbb{E}\left[
\left|Z_{t,\mathrm{ref}}^i
-\widetilde Z_{t,\Delta}^i\right|
\right]
\]
and
\[
e_2(\Delta)
=
\mathbb{E}\left[
\sup_{0\leq t\leq 1}
\left|Z_{t,\mathrm{ref}}^i
-\widetilde Z_{t,\Delta}^i\right|
\right].
\]

\begin{figure}[htbp]
    \centering
    \includegraphics[width=0.95\textwidth]
    {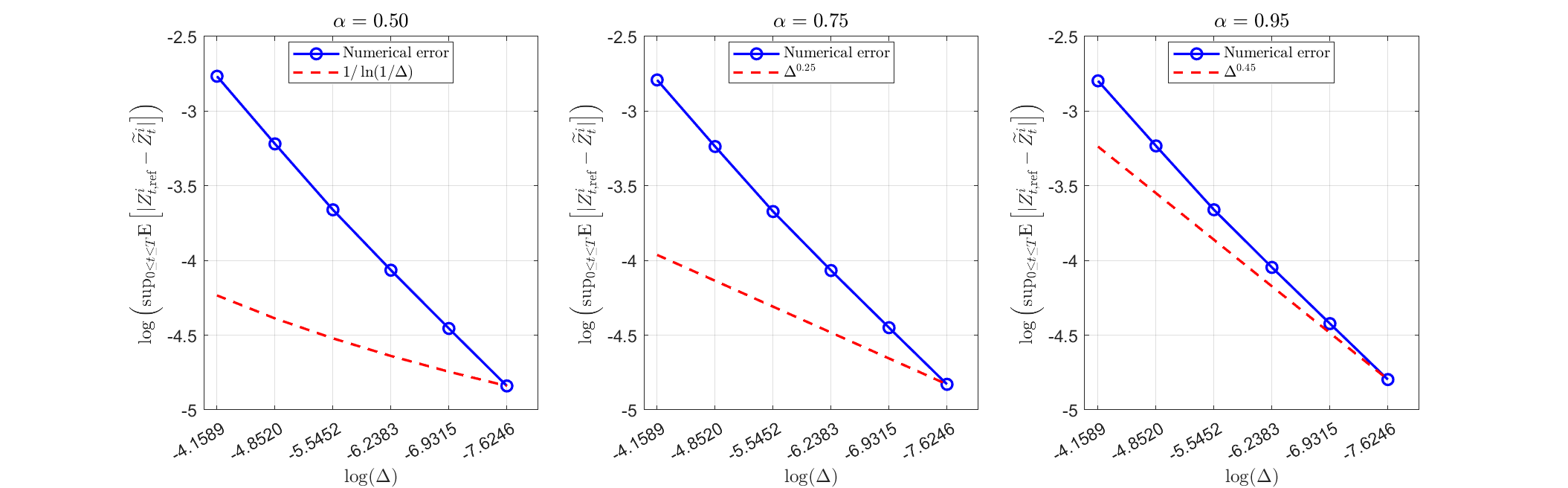}
    \caption{The logarithmic plot of the $L^1$ error $e_1(\Delta)$.}
    \label{fig:supE-error}
\end{figure}

\begin{figure}[htbp]
    \centering
    \includegraphics[width=0.95\textwidth]
    {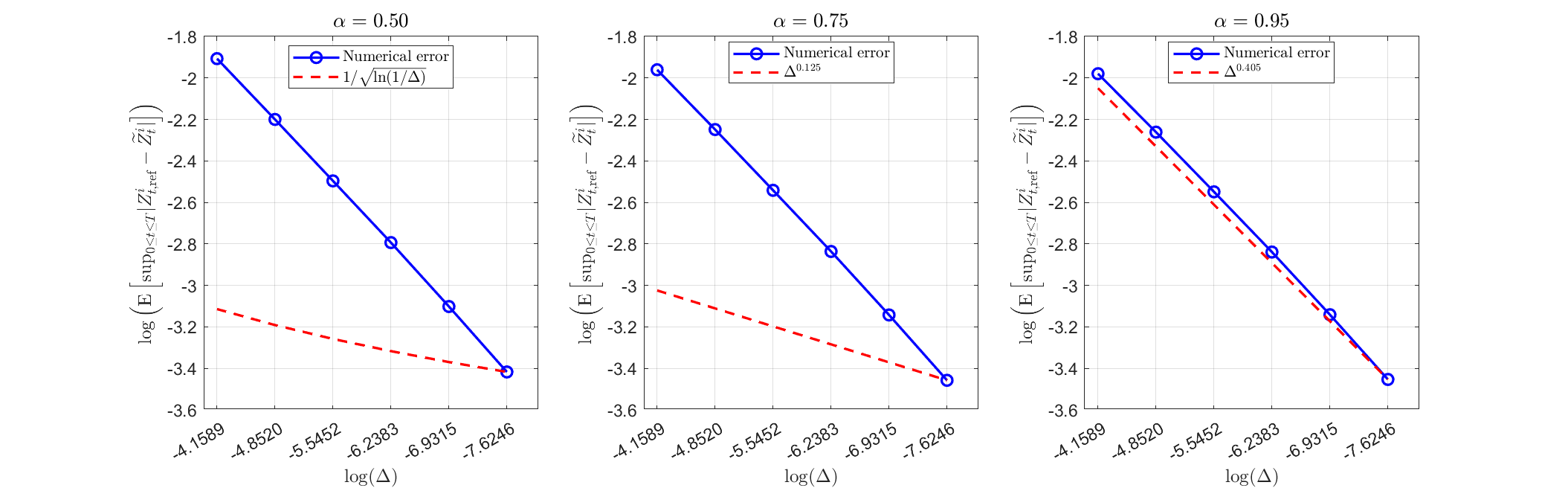}
    \caption{The logarithmic plot of the $L^1$ uniform error
    error $e_2(\Delta)$.}
    \label{fig:Esup-error}
\end{figure}

Figures~\ref{fig:supE-error} and~\ref{fig:Esup-error} show the
corresponding errors for the three H\"older exponents. The horizontal
and vertical axes are displayed on logarithmic scales. For
$\alpha=0.50$, the reference functions are
$1/\ln(1/\Delta)$ and $1/\sqrt{\ln(1/\Delta)}$ for the two error
quantities, respectively. For $\alpha=0.75$ and $\alpha=0.95$, the
reference orders are $\Delta^{0.25}$ and $\Delta^{0.45}$ for
$e_1(\Delta)$, and $\Delta^{0.125}$ and $\Delta^{0.405}$ for
$e_2(\Delta)$, respectively. These reference curves are included
to compare the observed behaviour with the theoretical convergence
rates.

\subsection{Particle-number error}

We next examine the effect of the particle number on the numerical
error. We use the same equation~\eqref{eq:convergence-test} and fix
the time step at $\Delta=2^{-9}$ and the terminal time at $T=1$.
The particle numbers are
\[
N=128,\ 256,\ 512,\ 1024,\ 2048,\ 4096.
\]
For each $N$, we compare the numerical solution of the full
$N$-particle system with two coupled subsystems, each containing
$N/2$ particles. The initial values and Brownian increments are
coupled so that the particles in the full system are matched with
the corresponding particles in the two half-size systems.

More precisely, let $Z_1^{i,N,r}$ denote the terminal value of the
$i$-th particle in the full system in the $r$-th replication, and
let $\widetilde Z_1^{i,N,r}$ denote the corresponding value obtained
from the two half-size systems. We estimate the particle-number
error by
\begin{equation}\label{eq:particle-error}
\widehat e_N
=
\frac{1}{UN}
\sum_{r=1}^{U}\sum_{i=1}^{N}
\left|Z_1^{i,N,r}
-\widetilde Z_1^{i,N,r}\right|,
\end{equation}
where $U=300$ is the number of independent replications.

\begin{figure}[htbp]
    \centering
    \includegraphics[width=0.95\textwidth]
    {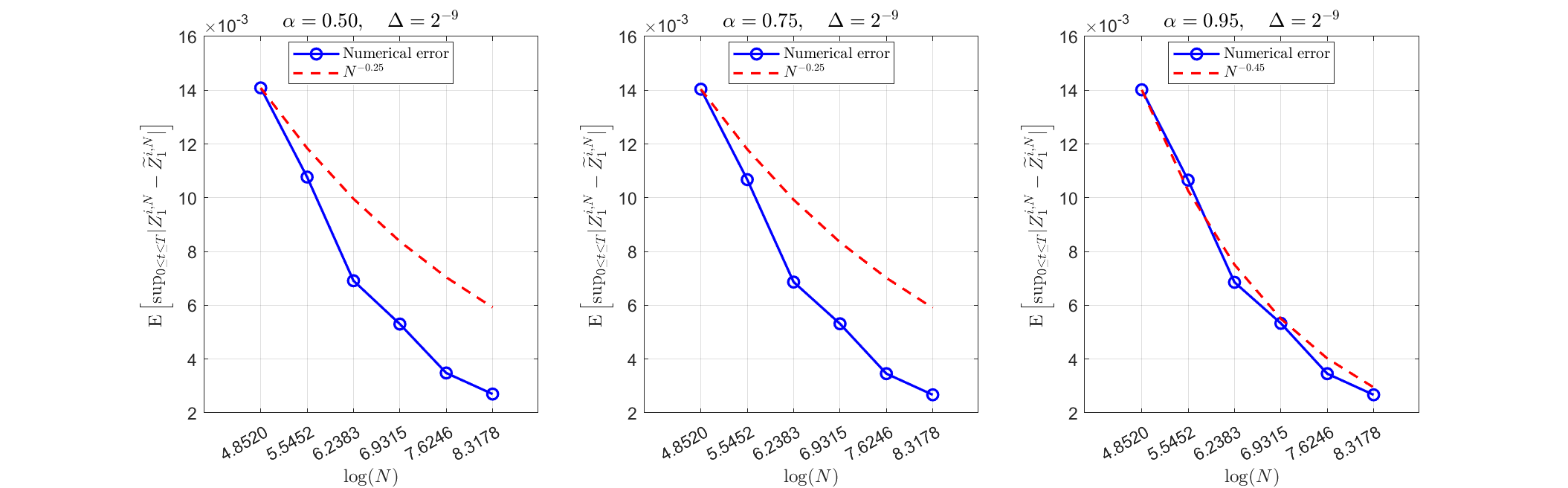}
    \caption{Strong convergence with respect to the number of particles}
    \label{fig:particle-error}
\end{figure}

Figure~\ref{fig:particle-error} presents the numerical errors for
the three H\"older exponents, together with the corresponding
reference curves. The reference orders are $N^{-0.25}$ for
$\alpha=0.50$ and $\alpha=0.75$, and $N^{-0.45}$ for
$\alpha=0.95$. The comparison provides a numerical indication of
the dependence of the particle-number error on the H\"older
regularity of the diffusion coefficient.

\subsection{Mean-square decay}

We next examine the long-time behavior of the splitting scheme
for the McKean--Vlasov SDE
\begin{equation}
\label{eq:stability-test}
\mathrm{d}X_t
=
\left(0.4\,\mathbb{E}[X_t]-X_t-X_t^3\right)\mathrm{d}t
+0.5|X_t|^\alpha\,\mathrm{d}W_t.
\end{equation}
As in the previous experiment, we consider
$\alpha=0.5, 0.75, 0.95$.

The drift is decomposed as
\[
f_1(x,\mu)
=
0.4\int_{\mathbb{R}}y\,\mu(\mathrm{d}y),
\qquad
f_2(x)=-x-x^3.
\]
The exact flow corresponding to the dissipative component is
\[
\Phi_\Delta(y)
=
\frac{ye^{-\Delta}}
{\sqrt{1+y^2(1-e^{-2\Delta})}}.
\]
For this experiment, we use the exponentially tamed diffusion
coefficient
\[
g_\Delta(t,x)
=
\frac{0.5|x|^\alpha}
{1+\sqrt{\Delta}\,e^{2\lambda t}
\left(1+0.5|x|^\alpha\right)},
\qquad \lambda=0.15.
\]

We take $N=10000$ particles and simulate the numerical solution
over the time interval $[0,16]$ with $\Delta=0.0125$. To
illustrate the decay of the second moment, we compute the empirical
quantity
\[
\frac{1}{N}\sum_{i=1}^{N}|Z_t^i|^2
\]
at each time step. In particular, we plot its logarithm divided
by time, namely
\[
\frac{1}{t}\log\left(
\frac{1}{N}\sum_{i=1}^{N}|Z_t^i|^2
\right),\qquad t>0.
\]

\begin{figure}[htbp]
    \centering
    \includegraphics[width=0.95\textwidth]
    {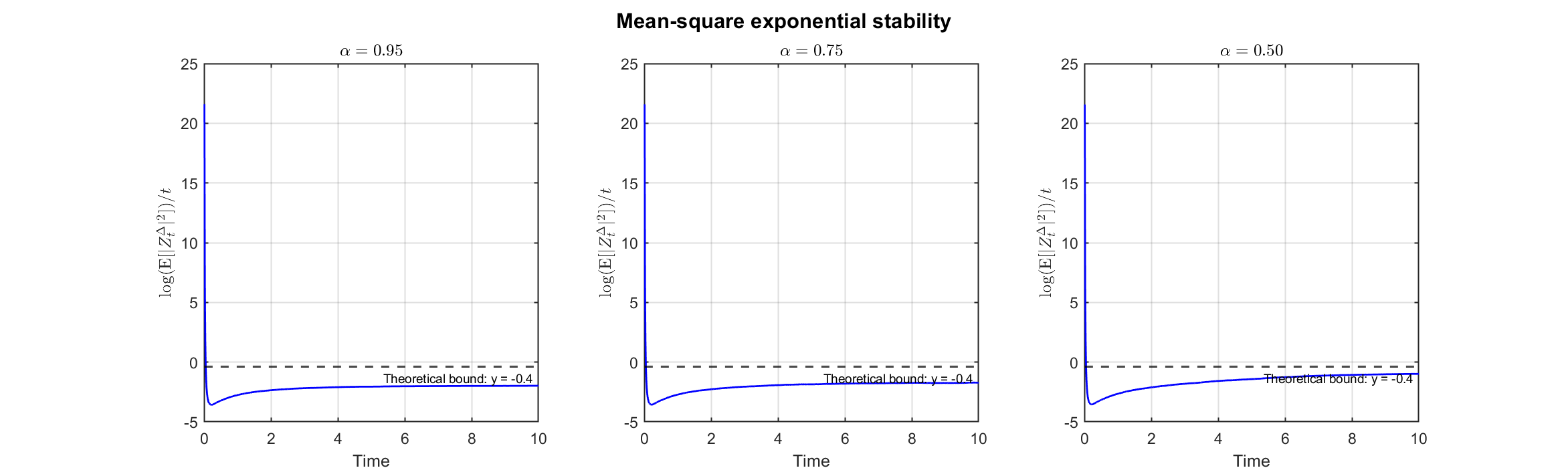}
    \caption{The logarithm of the empirical second moment divided
    by time for different H\"older exponents. The dashed line
    indicates the theoretical reference level $y=-0.4$.}
    \label{fig:mean-square-stability}
\end{figure}

For the given equation~\eqref{eq:stability-test}, the theoretical
exponential decay rate in the mean-square stability estimate is
$\eta_*=0.4$. Therefore, the horizontal reference line
$y=-0.4$ is included in the figure for all three H\"older
exponents. Figure~\ref{fig:mean-square-stability} displays the
resulting numerical curves and the theoretical reference line.

\begin{figure}[htbp]
    \centering
    \includegraphics[width=0.95\textwidth]
    {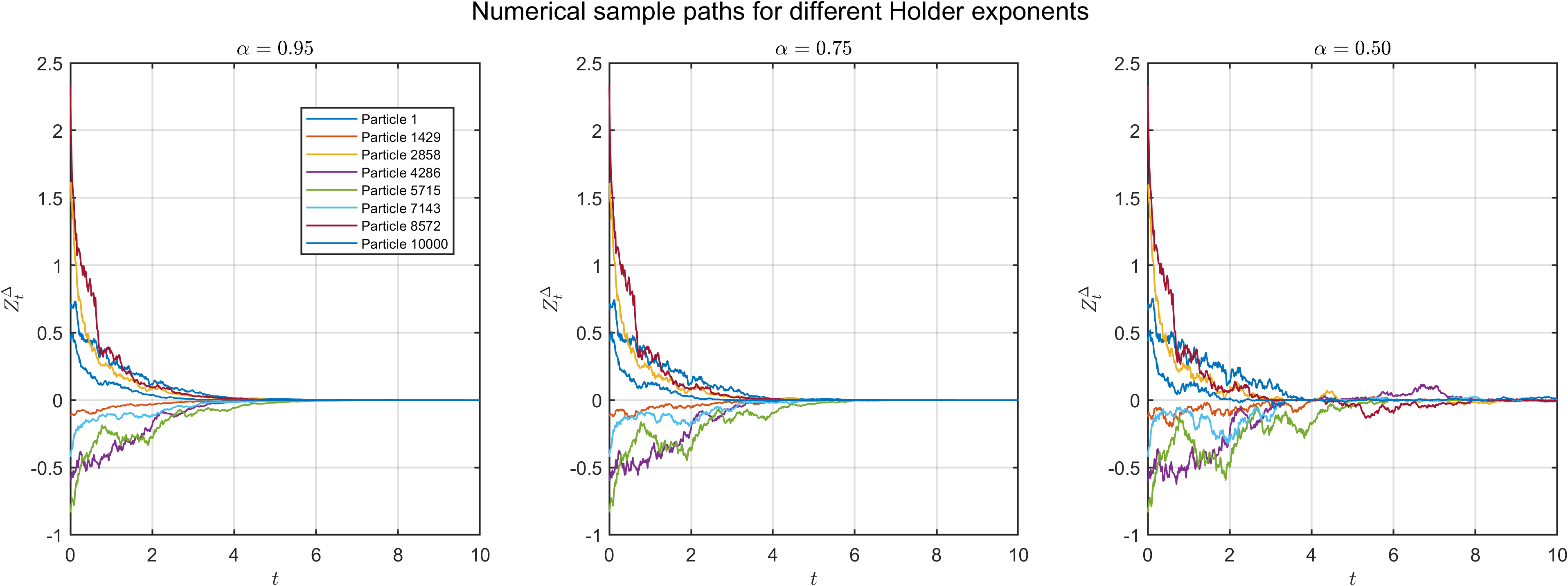}
    \caption{Numerical trajectories of eight selected
    particles among the $N=10000$ particles for different
    H\"older exponents.}
    \label{fig:squared-paths}
\end{figure}

\begin{figure}[htbp]
    \centering
    \includegraphics[width=0.95\textwidth]
    {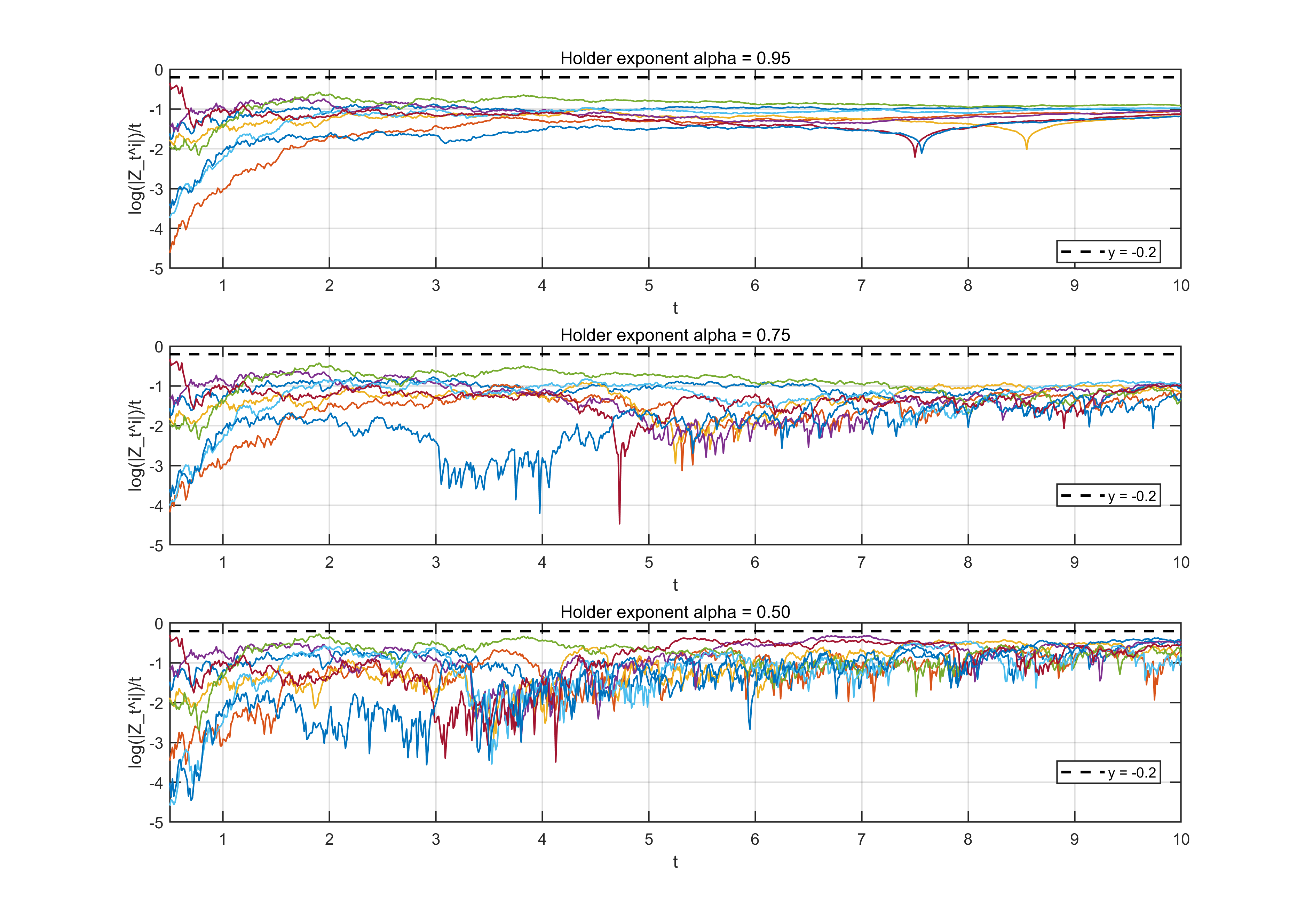}
    \caption{The logarithm of eight selected
    particles among the $N=10000$ particles divided
    by time for different H\"older exponents. The dashed line
    indicates the theoretical reference level $y=-0.2$.}
    \label{fig:LOG}
\end{figure}

To further illustrate the behaviour of individual numerical
trajectories, we select eight particles from the $N=10000$
particles and plot their values over time. Figure
\ref{fig:squared-paths} shows the trajectories
$Z_t^i$ for the three H\"older exponents while \ref{fig:LOG} indicates that $Z_t^i$ has almost surely exponential stability with rate $0.2$.

\newpage

\section{Concluding remark}\label{jiesu}

We obtain the strong convergence of a splitting Euler scheme for one-dimensional
McKean-Vlasov SDEs with a globally Lipschitz measure-dependent drift,
a superlinear dissipative drift component, and a H\"older continuous
diffusion coefficient. The exact flow of the dissipative component
allows us to obtain $L^1$ strong
convergence rate for
$0\le\beta\le\frac12$. $L^p$ uniform convergence over interval $[0,T]$ are also obtained with relatively lower rate. Under strict dissipative condition, both the exact solution and splitting Euler scheme with exponentially tamed diffusion are exponentially stable in $L^1$ and $L^2$, respectively.

\end{document}